\documentclass{birkjour}
\newtheorem{thm}{Theorem}[section]
\newtheorem{cor}[thm]{Corollary}
\newtheorem{lem}[thm]{Lemma}
\newtheorem{prop}[thm]{Proposition}
\theoremstyle{definition}
\newtheorem{defn}[thm]{Definition}
\theoremstyle{remark}
\newtheorem{rem}[thm]{Remark}
\newtheorem{ex}[thm]{Example}
\numberwithin{equation}{section}

\usepackage{enumerate}

\usepackage{mathrsfs}
\usepackage{amssymb}
\usepackage{amsmath}
\usepackage{stmaryrd}
\usepackage{lineno,hyperref} \modulolinenumbers[5]

\allowdisplaybreaks

\begin{document}
	
	%
	%
	%
	%
	%
	%
	%
	%
	%

	\title[Common fixed point theorems]
	{Common fixed point theorems for a commutative family of nonexpansive mappings in complete random normed modules}

	\author{Xiaohuan Mu}
	
	\address{%
		School of Mathematics and Statistics, \\
		Central South University,\\
		Changsha {\rm 410083}, \\
		China}
	
	\email{xiaohuanmu@163.com}
	
	
	\author{Qiang Tu}
	\address{%
		School of Mathematics and Statistics, \\
		Central South University,\\
		Changsha {\rm 410083}, \\
		China}
	
	\email{qiangtu126@126.com}
	
	\author{Tiexin Guo*}
	\address{%
		School of Mathematics and Statistics, \\
		Central South University,\\
		Changsha {\rm 410083}, \\
		China}
	
	\email{tiexinguo@csu.edu.cn}
	\thanks{*Corresponding author}
	
	\author{Hong-Kun Xu}
	\address{%
		School of Science, \\
		Hangzhou Dianzi University,\\
		 Hangzhou {\rm 310018}, \\
		 China\\
		College of Mathematics and Information Science,\\
		  Henan Normal University,\\
		Xinxiang, {\rm 453007},\\
		 China
	}

	\email{xuhk@hdu.edu.cn}

	
	
	\date{\today}
	
	\begin{abstract}
		In this paper, we first introduce and study the notion of random Chebyshev centers.
		Further, based on the recently developed theory of stable sets, we introduce
		the notion of random complete normal structure so that we can prove the two deeper theorems:
		one of which states that random complete normal structure is equivalent to random normal structure for 
		an $L^0$-convexly compact set in a complete random normed module; the other of which states that if $G$ is an $L^0$-convexly compact subset 
		with random normal structure of a complete random normed module, then every
		commutative family of nonexpansive mappings from $G$ to $G$ has a common fixed
		point.  We also consider the fixed point problems for isometric mappings in complete random normed modules.
		Finally, as applications of the fixed point theorems established in random normed modules, when the measurable selection theorems fail
		to work, we can still prove that a family of strong random  nonexpansive operators from 
		$(\Omega,\mathcal{F},P)\times C$ to $C$ has a common random fixed point, where $(\Omega,\mathcal{F},P)$ 
		is a probability space and $C$ is a weakly compact convex subset with normal structure of a Banach space.
	\end{abstract}
	
	\maketitle

	\setcounter{secnumdepth}{0}
	
	\section{Introduction}\label{intro}

	The fixed point theory of nonexpansive mappings in Banach spaces is one of the most attractive parts of
	metric fixed point theory since it not only has a strong interaction with geometry of Banach spaces but also
	is closely connected with the evolution equation governed by accretive operators, see the famous survey papers \cite{GR1984,Kirk1981,KS2001}.
	Let us first recall some known basic fixed point theorems concerning nonexpansive mappings, which is closely 
	related to our work in this paper. The well known Kirk fixed point theorem \cite{Kirk1965} states that if $K$ 
	is a nonempty weakly compact convex subset with normal
	structure of a Banach space, then every nonexpansive mapping from $K$ to itself
	has a fixed point.  In 1966, Belluce and Kirk \cite{BK1966} extended the
	Kirk fixed point theorem for any finite commutative family of
	nonexpansive self-mappings. Following \cite{BK1966}, in 1967, Belluce and Kirk \cite{BK1967} further introduced
	a strengthening notion of normal structure, called complete normal structure, and then proved that every
	commutative family of nonexpansive self-mappings acting on a nonempty weakly
	compact convex subset with complete normal structure of a Banach space has a
	common fixed point. In 1974, Lim \cite{lim1974} proved that complete normal
	structure is equivalent to normal structure for a weakly compact convex subset in
	locally convex spaces so that Lim \cite{Lim1974} can prove that
	the Kirk fixed point theorem holds true for any commutative family of
	nonexpansive self-mappings. In fact, in 1948, Brodskii and Milman \cite{BM1948} already proved that any family of 
	surjective isometric self-mappings acting on such a set $K$ as in Kirk fixed point theorem has a common fixed 
	point in the Chebyshev center of $K$. For a single isometric (not necessarily surjective) self-mapping,
	Lim et al. \cite{Lim2003} proved an interesting result, namely, every isometric self-mapping
	acting on such a set $K$ has a fixed point in the Chebyshev center of $K$.

	\par The central purpose of this paper is to extend these classical results in \cite{lim1974,Lim1974,BM1948,Lim2003} 
	from a Banach space to a complete random normed module. The crux of our work is to solve the problem of how to introduce
	a fruitful notion of random complete normal structure for an almost surely bounded closed $L^{0}$-convex subset of a 
	complete random normed module. The success of this paper considerably benefits from the recent advance in fixed point theory in random 
	functional analysis \cite{GZWG20,GWXYC24}. Random functional analysis is based on the idea of randomizing the traditional 
	space theory. Random normed modules, as a central framework of random functional analysis, are a random generalization of normed
	spaces. Over the past 30 years, random normed modules have been deeply developed and have played an essential role in the 
	development of the theory of conditional (dynamic) risk measures \cite{FKV2009,Guo10,GZZ14,Guo5} and  nonsmooth differential 
	geometry on metric measure spaces \cite{Gigli2018,GMT24,LP19,LPV24,CLP2024} (see \cite{Gigli2018,GMT24} for Gigli's independent contribution). 
	As is stated in \cite{GZWG20,GWXYC24}, one of the recent central tasks of random functional analysis is to extend 
	the classical fixed point theory in functional analysis to random functional analysis, which is not only of interest in its own right but also required to meet the needs of dynamic financial mathematics (for example, for the study of dynamic Nash equilibrium and dynamic optimal portfolio), and the challenge in 
	achieving this goal lies in overcoming noncompactness since the closed $L^{0}$-convex subsets frequently occurring in the theoretic developments
	and their financial applications are not compact in general. The classical Kirk fixed point theorem involves the two key notions --- normal structure 
	and weakly compact convex set. In the case of a complete random normed module, it is not very difficult to introduce the notion of random normal
	structure for a closed $L^{0}$-convex subset as in \cite{GZWG20}, but it is a delicate problem to speak of weak compactness for a closed
	$L^{0}$-convex subset. Since a random normed module is often endowed with the $(\varepsilon,\lambda)$-topology and is a metrizable but not locally
	convex topological module, it makes no sense to speak of weak compactness for a closed $L^{0}$-convex subset. Owing to \v{Z}itkovi\'{c}'s contribution 
	\cite{Zit2010}, he introduced the notion of convex compactness for a closed convex subset of a locally nonconvex space and developed some  powerful
	tools for the study of several important problems in finance and economics. It is motivated by the work in \cite{Zit2010} that Guo \cite{GZWW20}
	introduced the notion of $L^{0}$-convex compactness for a closed $L^{0}$-convex subset in a complete random normed module, and further established
	the characterization for a closed $L^{0}$-convex subset to be $L^{0}$-convexly compact by means of the theory of random conjugate spaces. Therefore,
	the notion of an $L^{0}$-convexly compact set is a proper substitution for the notion of a weakly compact convex subset of a Banach space, which
	had made Guo et al. smoothly generalize the Kirk fixed point theorem from a Banach space to a complete random normed module in \cite{GZWG20}. 
	Recently, based on the notion of $\sigma$-stable sets introduced in \cite{Guo10}, Guo et al. \cite{GWXYC24} developed the theory of random 
	sequentially compact sets and established a noncompact Schauder fixed point theorem, which directly leads to a more abstract development of 
	stable sets and  stably compact sets in \cite{GMT24,GWT23}. It is the development of the abstract stable sets that motivates us to introduce 
	the desired notion of random complete normal structure.
	
	\par The remainder of this paper is organized as follows: Section \ref{sec.1} provides some prerequisites and further presents the main results.  
	Section \ref{sec.2} is devoted to the proof of Theorem \ref{thm.1.14}, namely, proving the equivalence between random complete normal 
	structure and random normal structure for an $L^0$-convexly compact set in a complete random normed module. Section \ref{sec.3} is devoted to 
	the proof of Theorem \ref{thm.1.15}, which establishes a common fixed point theorem for a commutative family of  nonexpansive mappings 
	in a complete random normed module. Section \ref{sec.4} is devoted to the proofs of two fixed point Theorems \ref{thm.1.16} and \ref{thm.1.17} 
	for isometric mappings in a complete random normed module .  Section \ref{sec.5} is devoted to the applications of the main results to 
	random fixed point theorems for random nonexpansive operators.  Finally, Section \ref{sec.6} concludes with some remarks pointing out the future possible 
	works.

	\setcounter{secnumdepth}{3}
	
	\section{Prerequisites and main results}\label{sec.1}
	
	This section is divided into the following five subsections in order to clearly state some prerequisites and the main results
	of this paper.
	
	\subsection{$B_{\mathcal{F}}$-stable sets and consistent families with $B_{\mathcal{F}}$-stable index}\label{subsec.1.1}
	
	Throughout this paper, $\mathbb{K}$ denotes the scalar field $\mathbb{R}$ of real numbers or
	$\mathbb{C}$ of complex numbers, $(\Omega,\mathcal{F},P)$ a given probability space, $\mathbb{N}$ the set of positive integers, 
	$\mathbb{R}_{+}$ the set of nonnegative real numbers,
    $L^0(\mathcal{F},\mathbb{K})$ the algebra of equivalence classes of $\mathbb{K}$-valued $\mathcal{F}$-measurable random variables 
	on $(\Omega,\mathcal{F},P)$ (here, two random variables are said to be equivalent if they are equal almost surely), 
	$L^0(\mathcal{F}):=L^0(\mathcal{F},\mathbb{R})$ and $\bar{L}^0(\mathcal{F})$ the set of equivalence classes of extended
	real-valued $\mathcal{F}$-measurable random variables on $(\Omega,\mathcal{F},P)$. Besides, $I_A$ denotes the characteristic 
	function of $A$ for any $A \in	\mathcal{F}$ and $\tilde{I}_A$ denotes the equivalence class of $I_A$. 
	
	\par In the remainder of this paper, we always denote the set of countable partitions of $\Omega$ to $\mathcal{F}$ 
	by $\Pi_{\mathcal{F}}$, where a countable partition of $\Omega$ to $\mathcal{F}$ is a disjoint sequence $\{A_{n},n\in \mathbb{N}\}$ 
	in $\mathcal{F}$ such that $ \cup_{n=1}^{\infty} A_{n}=\Omega $.
	
	\par It is well known from \cite{DS58} that $\bar{L}^0(\mathcal{F})$ is a complete lattice under the partial order $\xi \leq \eta$
	iff $\xi^0(\omega) \leq \eta^0(\omega)$ for almost all $\omega \in \Omega$
	(briefly, $\xi^0(\omega) \leq \eta^0(\omega)$ a.s.), where $\xi^0$ and $\eta^0$
	are arbitrarily chosen representatives of $\xi$ and $\eta$ in
	$\bar{L}^0(\mathcal{F})$, respectively. Specially, the sublattice $(L^0(\mathcal{F}),\leq)$ is a Dedekind complete lattice. 
	For a nonempty subset $H$ of $\bar{L}^0(\mathcal{F})$, we usually use $ \bigvee H $ and $ \bigwedge H $ for the supremum and infimum of $H$, respectively. 
	
	\par  For any $A\in \mathcal{F}$, we
	always use the corresponding lowercase letter $a$ for the equivalence class $[A]$
	of $A$ (two elements $A$ and $D$ in $\mathcal{F}$ are said to be equivalent if
	$P(A\Delta D)=0$, where $A\Delta D=(A\setminus D)\cup(D\setminus A)$ stands
	for the symmetric difference of $A$ and $D$). For any $A,B\in \mathcal{F}$, define $a\wedge b=[A\cap B]$, $a\vee b=[A\cup B] $ and $a^{c}=[A^{c}]$, where $A^c$ denotes
	the complement of $A$. Then, it is well known from \cite{MKB89} that $B_{\mathcal{F}}=\{a=[A]:A\in \mathcal{F}\}$ is a complete Boolean algebra 
	with $1:=[\Omega]$ and $0:=[\emptyset]$, called the measure algebra associated with $(\Omega,\mathcal{F},P)$. A nonempty subset $\{a_i,i\in I\}$ of
	$B_{\mathcal{F}}$ is called a partition of unity if $\vee_{i\in I} a_i = 1$ and $a_i \wedge
	a_j = 0$ where $i\neq j$. It should be also noticed that $\{i\in I : a_i >
	0\}$ must be at most countable for any partition $\{a_i,i\in I\}$ of unity in
	$B_{\mathcal{F}}$. It is easy to check that,  $\{a_{n},n\in \mathbb{N}\}$ is a partition of unity in
	$B_{\mathcal{F}}$  iff  there exists $\{A_{n},n\in \mathbb{N}\}\in \Pi_{\mathcal{F}}$ such that $a_{n}=[A_{n}]$ for each $n\in \mathbb{N}$.
	
	\par For a random normed module $(E,\|\cdot\|)$ (see Subsection \ref{subsec.1.2} below for its definition and its $ (\varepsilon,\lambda) $-topology), 
	the $L^{0}$-norm can induce the two kinds of topologies --- the $ (\varepsilon,\lambda) $-topology and the locally $L^{0}$-convex topology, 
	the latter is much stronger than the former, it is in order to establish the inherent connection between the two kinds of topologies that 
	Guo \cite{Guo10} introduced the following notion of a $ \sigma $-stable set. First, let us recall the notion of a regular $L^{0}$-module. 
	A left module over the algebra $L^0(\mathcal{F},\mathbb{K})$ (briefly, an $L^0(\mathcal{F},\mathbb{K})$-module) is said to be regular if $E$ 
	has the following property: for any given two elements $x$ and $y$ in $E$, if there exists some $\{A_n,n\in \mathbb{N}\}\in \Pi_{\mathcal{F}}$ 
	such that $\tilde{I}_{A_{n}}x=\tilde{I}_{A_{n}}y$ for each $n\in \mathbb{N}$, then $x=y$. In the remainder of this paper, 
	we always assume that all the $L^0(\mathcal{F},\mathbb{K})$-modules occurring in this paper are regular, the assumption is not too restrictive since all 
	random normed modules are regular.

	\begin{defn}\label{defn.1.1}
		Let $E$ be an $L^0(\mathcal{F},\mathbb{K})$-module and $G$ be a nonempty
		subset of
		$E$. $G$ is said to be finitely stable if $\tilde{I}_{A}x+\tilde{I}_{A^{c}}y\in G$ for
		any $x,y\in G$ and any $A\in \mathcal{F}$. $G$ is said to be $\sigma$-stable
		(or to have the countable concatenation property) if for each sequence
		$\{x_n, n\in \mathbb{N} \}$ in $G$ and each $\{A_n,n\in
		\mathbb{N}\} \in \Pi_{\mathcal{F}}$, there exists some $x\in G$ such
		that $\tilde{I}_{A_n}x=\tilde{I}_{A_n}x_n$ for each $n\in \mathbb{N}$ ($x$ is
		unique since $E$ is assumed to be regular, usually denoted by
		$\sum_{n=1}^{\infty}\tilde{I}_{A_n}x_n$, called the countable concatenation of
		$\{x_n,n\in
		\mathbb{N}\}$ along $\{A_n,n\in \mathbb{N}\}$). By the way, if $G$ is
		$\sigma$-stable and  $H$ is a nonempty subset of $G$, then
		$\sigma(H):=\{\sum_{n=1}^{\infty}\tilde{I}_{A_n}h_n: \{h_n,n\in \mathbb{N}\}$
		is a sequence in $H$
		and $\{A_n,n\in \mathbb{N}\}\in \Pi_{\mathcal{F}}$\} is called the
		$\sigma$-stable hull of $H$.
	\end{defn}

	\par The notion of a $\sigma$-stable set depends on the structure of an $L^0(\mathcal{F},\mathbb{K})$-module, but the following  
	notion of a $B_{\mathcal{F}}$-stable set, as an abstract generalization of the notion of a $\sigma$-stable set, was introduced 
	in \cite{GWT23,GMT24} for the development of the theory of stably compact sets, which will be needed in order to introduce the 
	notion of random complete normal structure in Subsection \ref{subsec.1.2} of this paper.

	\begin{defn}\label{defn.1.2}
		Let $X$ be a nonempty set. An equivalence relation $\sim$ on
		$X\times B_{\mathcal{F}}$ (denote the equivalence class of $(x,a)$ by $x|a$ for
		each $(x,a)\in X\times B_{\mathcal{F}}$) is said to be regular if the following
		conditions are satisfied:
		\begin{itemize}
			\item [(1)] $x|a=y|b$ implies $a=b$;
			\item [(2)]	$x|a=y|a$ implies $x|b=y|b$ for any $b\leq a$;
			\item [(3)]	$\{a\in B_{\mathcal{F}}: x|a=y|a\}$ has a greatest element for
			any $x$ and $y$ in $X$;
			\item [(4)]	$x|1=y|1$ implies $x=y$.
		\end{itemize}
		In addition, given a regular equivalence relation on $X\times
		B_{\mathcal{F}}$, a nonempty subset $G$ of $X$ is said to be
		$B_{\mathcal{F}}$-stable with respect to $\sim$ if, for each sequence
		$\{x_n,n\in \mathbb{N}\}$ in $G$ and each partition $\{a_n,n\in \mathbb{N}\}$ of
		unity in $B_{\mathcal{F}}$, there exists $x\in G$ such that $x|a_n=x_n|a_n$ for
		each $n\in \mathbb{N}$ (by (3) and (4) it is easy to see that such an $x$ is
		unique, denoted by $\sum_{n=1}^{\infty}x_{n}|a_{n}$).  By the way, if $G$ is
		$B_{\mathcal{F}}$-stable and $H$ is a nonempty subset of $G$, then
		$B_{\sigma}(H):=\{\sum_{n=1}^{\infty}h_n|a_{n}: \{h_n,n\in \mathbb{N}\}$
		is a sequence in $H$
		and $\{a_n,n\in \mathbb{N}\}$ is a partition of unity in $B_{\mathcal{F}}$\}
		is called the
		$B_{\mathcal{F}}$-stable hull of $H$ with respect to $\sim$.
	\end{defn}
	
	\begin{ex}\label{ex.1.3}
		Let $E$ be an $L^0(\mathcal{F},\mathbb{K})$-module. Define an
		equivalence relation $\sim$ on $E\times B_{\mathcal{F}}$ by $(x,a)\sim (y,b)$
		iff $a=b$ and $\tilde{I}_{A}x=\tilde{I}_{A}y$, then  $\sim$ is regular. It is easy to check that a nonempty subset $G$ is $ \sigma $-stable iff
		$ G $ is  $B_{\mathcal{F}}$-stable, and in this case, for any sequence $\{x_n,n\in \mathbb{N}\}$ in $G$ and any partition $\{a_n,n\in
		\mathbb{N}\}$ of unity in $B_{\mathcal{F}}$,
		$\sum_{n=1}^{\infty}x_{n}|a_{n}=\sum_{n=1}^{\infty}\tilde{I}_{A_{n}}x_{n}$. Consequently, the notion of a $ \sigma $-stable set is a special case of that of a $B_{\mathcal{F}}$-stable set.
	\end{ex}
	
	\par Companying Definition \ref{defn.1.2}, Definition \ref{defn.1.4} below was also introduced in \cite{GWT23}.
	
	\begin{defn}\label{defn.1.4}
		Let $E$ be an $L^{0}(\mathcal{F},\mathbb{K})$-module and $G$ be a
		$\sigma$-stable subset of $E$. For any sequence of
		nonempty subsets $\{G_n, n \in \mathbb{N}\}$ of $G$ and any $\{A_n, n \in
		\mathbb{N}\}\in \Pi_{\mathcal{F}}$, $\sum^{\infty}_{n=1} \tilde{I}_{A_n} G_n:$
		$=\{\sum^{\infty}_{n=1} \tilde{I}_{A_n} x_n: x_n\in G_n, \forall~ n\in
		\mathbb{N}\}$ is called the countable concatenation of $\{G_n,n\in\mathbb{N}\}$
		along $\{A_n,n\in \mathbb{N}\}$.
		For a nonempty family $\mathcal{E}$ of $\sigma$-stable subsets of $G$,
		$\sigma(\mathcal{E}):=\{\sum_{n=1}^{\infty}\tilde{I}_{A_n}G_{n}:
		\{G_{n},~n\in\mathbb{N}\}$~is a sequence in ~$\mathcal{E}$ and
		$\{A_{n},~n\in\mathbb{N}\}\in \Pi_{\mathcal{F}}\}$ is called the
		$\sigma$-stable hull
		of $\mathcal{E}$; if $\sigma(\mathcal{E})=\mathcal{E}$, then $\mathcal{E}$ is
		said to be $\sigma$-stable.
	\end{defn}
	
	\begin{rem}\label{rem.1.5}
		Let $ E $, $ G $ and $\mathcal{E}$ be the same as in Definition \ref{defn.1.4}. Define the equivalence relation $\sim$ on $\mathcal{E}\times B_{\mathcal{F}}$ by $(G_{1},a)\sim (G_{2},b)$
		iff $a=b$ and $\tilde{I}_{A}G_{1}=\tilde{I}_{A}G_{2}$, then it is easy to check that $\sim$ is regular. Further, it is also easy to verify that 
		$\mathcal{E}$ is $\sigma$-stable  iff  $\mathcal{E}$ is $B_{\mathcal{F}}$-stable with respect to $\sim$, and in this case, 
		$ \sum_{n=1}^{\infty}G_{n}|a_{n}=\sum_{n=1}^{\infty}\tilde{I}_{A_{n}}G_{n}$ for any sequence $\{G_n,n\in \mathbb{N}\}$ in $\mathcal{E}$ 
		and any partition $\{a_n,n\in\mathbb{N}\}$ of unity in $B_{\mathcal{F}}$ (please bear in mind $a_{n}=[A_{n}]$ according to our convention). 
	\end{rem}
	
	\par Let us conclude Subsection \ref{subsec.1.1} with the following notion of a consistent family with $B_{\mathcal{F}}$-stable index.
	
	\begin{defn}\label{defn.1.6}
		Let $E$ be an $L^{0}(\mathcal{F},\mathbb{K})$-module, $G$ a $\sigma$-stable subset of $E$ and $\Lambda$ a $B_{\mathcal{F}}$-stable set. 
		A nonempty family $\mathcal{E}:= \{G_{\alpha},\alpha\in \Lambda\} $ of $\sigma$-stable subsets of $G$ is called a consistent family with 
		$B_{\mathcal{F}}$-stable index if $ G_{\sum_{n=1}^{\infty}\alpha_{n}|a_{n}}=\sum_{n=1}^{\infty}\tilde{I}_{A_{n}}G_{\alpha_{n}} $ for any 
		sequence $\{G_{\alpha_{n}},n\in \mathbb{N}\}$ in $\mathcal{E}$ and any partition $\{a_n,n\in\mathbb{N}\}$ of unity in $B_{\mathcal{F}}$.  
		In addition, if $ \Lambda $ is a $B_{\mathcal{F}}$-stable directed set, then the consistent family $ \{G_{\alpha},\alpha\in \Lambda\} $ is 
		called a consistent net. Finally,  a consistent net $ \{G_{\alpha},\alpha\in \Lambda\} $ is said to be decreasing if $ G_{\alpha}\subset G_{\beta} $ 
		when $ \alpha\geq \beta $.
	\end{defn}
	
	\par By Remark \ref{rem.1.5}, a consistent family with $B_{\mathcal{F}}$-stable index must be $\sigma$-stable in the sense of Definition \ref{defn.1.4} 
	since  $\sum_{n=1}^{\infty}\tilde{I}_{A_{n}}G_{\alpha_{n}}=\sum_{n=1}^{\infty}G_{\alpha_{n}}|a_{n} $ is just $G_{\sum_{n=1}^{\infty}\alpha_{n}|a_{n}}\in \mathcal{E}$.

	\subsection{Random complete normal structure and random normal structure}\label{subsec.1.2}
	
	The following notion of a random normed module was introduced by Guo in \cite{Guo92,Guo93}, see also \cite{Guo10}. 
	The notion of an $L^0$-normed $L^0$-module independently introduced by Gigli \cite{Gigli2018} in the context of  
	nonsmooth differential geometry on metric measure spaces is equivalent to the notion of a random normed module. 
	Throughout this paper, $L^0_{+}(\mathcal{F})=\{\xi\in L^0(\mathcal{F}): \xi\geq 0\}$.

	\begin{defn}\label{definition2.2}
		An ordered pair $(E,\| \cdot \|)$ is called a random normed module $($briefly,
		an $RN$ module$)$ over $\mathbb{K}$ with base
		$(\Omega,\mathcal{F},P)$ if $E$ is  an
		$L^0(\mathcal{F},\mathbb{K})$--module and $\| \cdot \|$ is a mapping from $E$
		to $L^0_+(\mathcal{F})$ such that the following conditions are
		satisfied:
		\begin{enumerate}[(1)]
			\item $\| \xi x \| = |\xi| \|x\|$ for any $\xi \in L^0(\mathcal{F},
			\mathbb{K})$ and any $x \in E$;
			\item $\|x+y\| \leq \|x\| + \|y\|$ for any $x$ and $y$ in $E$;
			\item $\|x\| = 0$ implies $x = \theta$ $($the null element in $E)$.
		\end{enumerate}
		As usual, $\| \cdot \|$ is called the $L^0$-norm on $E$. If $\| \cdot \| :
		E \to L^0_+(\mathcal{F})$ only satisfies (1) and (2), then it is called an
		$L^0$-seminorm on $E$.
	\end{defn}

	\par As mentioned in Subsection \ref{subsec.1.1},  for an $RN$ module $(E,\|\cdot\|)$, 
	the $L^{0}$-norm $\|\cdot\|$ can induce the $ (\varepsilon,\lambda) $-topology as follows.

	\begin{defn}\label{definition2.3}
		Let $(E,\|\cdot\|)$ be an $RN$ module over $\mathbb{K}$ with base
		$(\Omega,\mathcal{F},P)$. For any real numbers $\varepsilon$ and
		$\lambda$ with $\varepsilon>0$ and $0< \lambda <1$, let
		$N_{\theta}(\varepsilon, \lambda)=\{x\in E: P\{\omega\in \Omega:
		\|x\|(\omega)<\varepsilon\}>1-\lambda\}$, then
		$\mathcal{U}_{\theta}:=\{N_{\theta}(\varepsilon, \lambda): \varepsilon>0,
		0<\lambda <1\}$ forms a local base of some metrizable Hausdorff
		linear topology for $E$, called the $(\varepsilon, \lambda)$-topology, denoted
		by $\mathcal{T}_{\varepsilon, \lambda}$.
	\end{defn}
	
	\par The $(\varepsilon,\lambda)$-topology is an abstract generalization of the
	topology of convergence in probability, in fact,
	a sequence $\{x_n, n \in \mathbb{N}\}$ in an $RN$ module converges in the
	$(\varepsilon,\lambda)$-topology to $x$ if and only if $\{\|x_n -x\|,n \in
	\mathbb{N}\}$ converges in probability to 0.

	\par When $(\Omega,\mathcal{F},P)$ is trivial, namely
	$\mathcal{F}=\{\Omega,\emptyset\}$, an $RN$ module $(E,\|\cdot\|)$ over $\mathbb{K}$ with base
	$(\Omega,\mathcal{F},P)$ reduces to an ordinary
	normed space over $\mathbb{K}$, and the $(\varepsilon,\lambda)$-topology to the usual norm topology.
	The simplest nontrivial $RN$ module is
	$(L^{0}(\mathcal{F},\mathbb{K}),|\cdot|)$,
	where $|\cdot|$ is the absolute value mapping.

	\par To extend the results of \cite{lim1974,Lim1974,BM1948,Lim2003} from Banach spaces
	to complete $RN$ modules, the main challenge is to introduce the notion of random complete 
	normal structure in $RN$ modules.  Before doing so, let us first recall the notions of normal 
	structure \cite{BM1948} and complete normal structure \cite{BK1967} in Banach spaces.
	
	Let $(B,\|\cdot\|)$ be a Banach sapce, $F$ a nonempty subset of $B$ and $K$ a  nonempty
	bounded subset of $B$. Define $r(K,\cdot):B\rightarrow \mathbb{R}_{+}$ by
	$$r(K,x)=\sup\{\|x-y\|:y\in H\},~\forall x\in B.$$
	Then $r(K,F)=\inf\{r(K,x): x\in F\}$
	and
	$c(K,F)=\{x\in F: r(K,x)=r(K,F)\}$
	are called the Chebyshev radius
	and  Chebyshev center of $K$ with respect to $F$, respectively, which were introduced in \cite{BM1948}.
	In particular, $r(K,K)$ and $c(K,K)$ are called the Chebyshev
	radius and Chebyshev center of $K$, respectively. Further, a nonempty bounded closed convex set $K$ in 
	a Banach space $(B,\|\cdot\|)$ is said to have normal structure if for every nonempty closed convex subset 
	$F$ of $K$, with $ F $ containing more than one point, has a point $x\in F$ such that $r(F,x)<\delta(F)$, 
	where $\delta(F)=\sup\{\|z-y\|:y,z\in F\}$ is the diameter of $ F $.
	
	\par The following notion of complete normal structure was introduced in \cite{BK1967} and once played 
	an important role in the study of common fixed point theorems.

	\begin{defn}\label{defn.cns}
		Let $(B,\| \cdot \|)$ be a Banach space. A bounded closed convex subset $K$ of
		$B$ is said to have complete normal structure if every closed convex subset $W$
		of $K$, with $ W $ containing more than one point, satisfies the following condition:
		\begin{enumerate}
			\item [$(\ast)$] For every decreasing net $\{W_\alpha, \alpha \in \Lambda\}$
			of nonempty subsets of $W$, if $r(W_\alpha,W)=r(W,W) ~\forall~ \alpha \in \Lambda$, then 
			the closure of $\bigcup_{\alpha \in \Lambda} c(W_\alpha,W) $ is a nonempty
			proper subset of $W$.
		\end{enumerate}
	\end{defn}

	\par Now, we will introduce the notion of random Chebyshev centers in $RN$ modules as
	follows. Let us first recall that a subset $G$ of an $RN$ module $(E,\|\cdot\|)$ with base
	$(\Omega,\mathcal{F},P)$ is said to almost surely (briefly, a.s.) bounded if $\bigvee\{\|g\|:g\in G\}\in
	L^0_{+}({\mathcal{F}})$.

	\begin{defn}\label{defn. random Chebyshev center}
		Let $(E,\|\cdot\|)$ be an $RN$ module over $\mathbb{K}$ with base
		$(\Omega,\mathcal{F},P)$, $G$ a nonempty subset of $E$ and $H$ an a.s.
		bounded nonempty subset of $E$. Define $R(H,\cdot):E\rightarrow L^0_{+}(\mathcal{F})$ by
		$$R(H,x)=\bigvee\{\|x-y\|:y\in H\},~\forall x\in E.$$
		Then $R(H,G)=\bigwedge\{R(H,x): x\in G\}$
		and
		$C(H,G)=\{x\in G: R(H,x)=R(H,G)\}$
		are called the random Chebyshev radius
		and random Chebyshev center of $H$ with respect to $G$, respectively.
		In particular, we call $R(H,H)$ and $C(H,H)$ the random Chebyshev
		radius and random Chebyshev center of $H$, respectively.
	\end{defn}

	\par Let $E$ be an $L^{0}(\mathcal{F},\mathbb{K})$-module and $G$
	be a nonempty subsets of $E$. $G$ is
	said to be $L^0$-convex if $\xi x+(1-\xi)y \in G$ for any $x,y\in G$ and any
	$\xi \in L^0_{+}(\mathcal{F})$ with $0\leq \xi \leq 1$.  Similarly, one can have the notion of the $L^{0}$-convex hull of $G$, denoted by
	$Conv_{L^{0}}(G)$.
	
	\par Let $\xi$ and $\eta$ be in $\bar{L}^0(\mathcal{F})$. As usual, $\xi >
	\eta$ means  $\xi\geq\eta$  but $\xi \neq \eta$. Besides, for any $A\in \mathcal{F}$ with $P(A)>0$, $\xi >
	\eta$ on $A$ means $\xi^0(\omega)>\eta^0(\omega)$ for almost all $\omega \in A$,
	where $\xi^0$ and $\eta^0$ are arbitrarily chosen representatives of $\xi$ and
	$\eta$, respectively. We always use
	$(\xi> \eta)$ for the set $\{\omega\in \Omega: \xi^{0}(\omega)>
	\eta^{0}(\omega)\}$ for any two representatives $\xi^0$ and $\eta^0$ of
	$\xi$ and $\eta$, respectively.  Although $(\xi> \eta)$ depends on the particular
	choice of $\xi^{0}$ and $\eta^{0}$, a careful reader can check that the subsequent
	notions and propositions involving $(\xi> \eta)$ are independent of the particular
	choice of $\xi^{0}$ and $\eta^{0}$.
	
	\par The following notion of random normal structure introduced by Guo in \cite{GZWG20} 
	was used to establish the random Kirk fixed point theorem.
	
	\begin{defn}\label{defn.rns}
		Let $(E,\|\cdot\|)$ be a $\mathcal{T}_{\varepsilon,\lambda}$-complete $RN$ module over $\mathbb{K}$ with base
		$(\Omega,\mathcal{F},P)$ and $G$ be a nonempty
		$\mathcal{T}_{\varepsilon,\lambda}$-closed $L^0$-convex subset of $E$. $G$ is
		said to have random normal structure if for each a.s. bounded
		$\mathcal{T}_{\varepsilon,\lambda}$-closed $L^0$-convex subset $H$ of $G$ such
		that $D(H):= \bigvee \{\|x-y\|: x,y\in H\}> 0$, there exists a nondiametral
		point $z\in H$, namely $R(H,z)< D(H)$ on $(D(H)> 0)$, where $D(H)$ is called the
		random diameter of $H$.
	\end{defn}
	
	\par One can observe that it is easy to generalize the notion of normal structure to that of 
	random normal structure, but it is completely another matter to generalize the notion of complete 
	normal structure to the following notion of random complete normal structure. Throughout this paper, for a subset 
	$G$ of an $RN$ module, $G^{-}_{\varepsilon,\lambda}$ denotes the closure of $G$ under the $(\varepsilon,\lambda)$-topology.
	
	\begin{defn} \label{defn.1.12}
		Let $( E , \| \cdot \|)$ be a $\mathcal{T}_{\varepsilon,\lambda}$-complete
		$RN$ module over $\mathbb{K}$ with base $(\Omega,\mathcal{F},P)$. An
		a.s. bounded $\mathcal{T}_{\varepsilon,\lambda}$-closed $L^0$-convex subset $G$
		of $E$ is said to have random complete normal structure if every
		$\mathcal{T}_{\varepsilon,\lambda}$-closed $L^0$-convex subset $W$ of $G$, with $ W $ 
		containing more than one point, satisfies the following condition: 
		\begin{enumerate}
			\item [$(\ast)$] For every decreasing consistent net $\{W_\alpha, \alpha \in
			\Lambda\}$ of $\sigma$-stable subsets of $W$, if $R(W_\alpha,W)=R(W,W) ~\forall~
			\alpha \in \Lambda$, then $[\bigcup_{\alpha \in \Lambda} C(W_\alpha,W)
			]^{-}_{\varepsilon,\lambda}\neq\emptyset$ and for any $B\subset(D(W)>0)$ with
			$P(B)>0$,
			there exists $y_{B}\in W$, such that $\tilde{I}_{B}y_{B} \notin \tilde{I}_{B}
			[\bigcup_{\alpha \in \Lambda} C(W_\alpha,W)
			]^{-}_{\varepsilon,\lambda} $.
		\end{enumerate}
	\end{defn}
	
	\par When $\mathcal{F}=\{\emptyset,\Omega\}$, Definition \ref{defn.1.12}, of course, 
	reduces to Definition \ref{defn.cns}.  The significant difference between Definition \ref{defn.1.12} 
	and Definition \ref{defn.cns} lies in that Definition \ref{defn.1.12} requires that each $W_\alpha$ 
	be $\sigma$-stable and $\{W_\alpha, \alpha \in\Lambda\}$ a  decreasing consistent net. Besides, 
	in Definition \ref{defn.1.12}, we are forced to consider the case for each $B\subset(D(W)>0)$ with $P(B)>0$. 
	Such a difference makes the proofs of Theorems \ref{thm.1.14} and \ref{thm.1.15} below more 
	involved than those of their prototypes \cite{lim1974,Lim1974}, and in particular we are required to frequently
	construct a  decreasing consistent net for the proof of Theorem \ref{thm.1.15}.

	\par Since a complete $RN$ module is generally a locally nonconvex
	space, so it does not make sense to speak of weak compactness for a closed
	$L^{0}$-convex subset. As a proper substitution for ordinary weak
	compactness, the notion of $L^{0}$-convex compactness was introduced and studied
	by Guo et al. in \cite{GZWW20}. The initial aim of \cite{GZWW20} is to generalize 
	the notion of convex compactness for a closed convex set in a linear topological 
	space to the notion of $L^{0}$-convex compactness for a closed $L^{0}$-convex subset 
	in a topological module over the topological algebra $ L^{0}(\mathcal{F},\mathbb{K}) $, 
	it turns out that $L^{0}$-convex compactness and convex compactness are equivalent for 
	a closed $L^{0}$-convex subset in a complete $RN$ module or a more general complete random 
	locally convex module, see \cite{GZWW20,WZZ22} for details. But playing a truly crucial role 
	in random functional analysis are closed $L^{0}$-convex subsets rather than generic closed 
	convex subsets, and the theory of random conjugate spaces helps establish the Jame's type 
	theorem characterizing $L^{0}$-convex compactness for closed $L^{0}$-convex subsets in \cite{GZWW20,WZZ22}. 
	Therefore, we would like to retain the terminology of $L^{0}$-convex compactness for closed $L^{0}$-convex subsets.
	
	\begin{defn}\label{defn.1.13}
		Let $(E,\|\cdot\|)$ be a $\mathcal{T}_{\varepsilon,\lambda}$-complete $RN$ module over $\mathbb{K}$ with base
		$(\Omega,\mathcal{F},P)$. A nonempty $\mathcal{T}_{\varepsilon, \lambda}$-closed
		$L^0$-convex subset $G$ of $E$ is said to be $L^0$-convexly compact (or, to have
		$L^0$-convex compactness) if every family of $\mathcal{T}_{\varepsilon,
			\lambda}$-closed $L^0$-convex subsets of $G$ has a nonempty intersection
		whenever the family has the finite intersection property.
	\end{defn}

	\par The first main result of this paper is Theorem \ref{thm.1.14} below, which can be 
	regarded as a proper generalization of \cite[Corollary 1]{lim1974} in the case of a Banach space.

	\begin{thm}\label{thm.1.14}
		Let $( E , \| \cdot \|)$ be a $\mathcal{T}_{\varepsilon,\lambda}$-complete
		$RN$ module over $\mathbb{K}$ with base $(\Omega,\mathcal{F},P)$ and
		$G$ be an $L^0$-convexly compact subset of $E$. Then $G$ has random complete
		normal structure iff $G$ has random normal structure.
	\end{thm}

	\subsection{A common fixed point theorem for a commutative family of nonexpansive mappings in a complete random normed module}\label{subsec.1.3}
	Let $(E,\|\cdot\|)$ be an $RN$ module, $G$ and $F$ be two nonempty subsets of
	$E$. The mapping $T: G\rightarrow F$ is said to be nonexpansive if  $\|Tx-Ty\|\leq \|x-y\| $ for any $x,y\in G$.

	\par The second main result of this paper is Theorem \ref{thm.1.15} below, which generalizes \cite[Theorem
	2.1]{BK1967} and \cite[Theorem 1]{Lim1974} from Banach spaces to $\mathcal{T}_{\varepsilon,\lambda}$-complete $RN$ modules. It 
	also extends the random Kirk fixed point theorem \cite{GZWG20} to the case of a commutative family of nonexpansive mappings.

	\begin{thm}\label{thm.1.15}
		Let $( E , \| \cdot \|)$ be a $\mathcal{T}_{\varepsilon,\lambda}$-complete
		$RN$ module over $\mathbb{K}$ with base $(\Omega,\mathcal{F},P)$, $G$
		an $L^0$-convexly compact subset with random normal structure of $E$ and
		$\mathcal{T}$ a commutative family of nonexpansive mappings from $G$ to $G$.
		Then $\mathcal{T}$ has a common fixed point.
	\end{thm}

	\subsection{A common fixed point theorem for a family of surjective isometric mappings in a complete random normed module}\label{subsec.1.4}

	As in the spirit of \cite{BM1948} and \cite{Lim2003}, namely, when the nonexpansive mappings are isometric, 
	we can locate fixed points in the Chebyshev center. Theorem \ref{thm.1.16} below is a random generalization  
	of \cite[Theorem 2]{Lim2003}. In fact, it is also a strengthening form of \cite[Theorem 3.8]{GZWG20} in the special 
	case when the mapping $T$ in \cite[Theorem 3.8]{GZWG20} is isometric.
	
	\begin{thm} \label{thm.1.16}
		Let $( E , \| \cdot \|)$ be a $\mathcal{T}_{\varepsilon,\lambda}$-complete
		$RN$ module over $\mathbb{K}$ with base $(\Omega,\mathcal{F},P)$ and
		$G$ be an $L^0$-convexly compact subset with random normal structure of $E$.
		Then every isometric mapping $T: G\rightarrow G$ (namely, $\|Tx-Ty\|= \|x-y\|$ for any $x,y\in G$) has a fixed point in
		$C(G,G)$.
	\end{thm}

	Theorem \ref{thm.1.17} below is a random generalization of
	Brodskii and Milman's result in \cite{BM1948}.
	
	\begin{thm} \label{thm.1.17}
		Let $( E , \| \cdot \|)$ be a $\mathcal{T}_{\varepsilon,\lambda}$-complete
		$RN$ module over $\mathbb{K}$ with base $(\Omega,\mathcal{F},P)$, $G$
		an $L^0$-convexly compact subset with random normal structure of $E$ and $ \mathcal{T} $ a family of surjective 
		isometric mappings from $G$ to $G$. Then $\mathcal{T}$ has a common fixed point  in $C(G,G)$.
	\end{thm}

	\subsection{Some applications to random fixed point theorems for random nonexpansive operators}\label{subsec.1.5}
	
	Random fixed point theory of random operators is aimed at giving the random versions of 
	classical fixed point theorems, which is required for the study of various classes of 
	random equations, see \cite{Han57,Bha1976} for a historical survey of random operator theory. 
	It already obtained an important advance in 1970s when Bharucha-Reid et al. 
	established the random  version of the Schauder fixed point theorem \cite{Bha1976}. Subsequently, random fixed point theory 
	further ushered in an explosive development \cite{Pa74,En78,En1978,Reich1978,Iton1979,Lin1988,Papa1990,Xu1990,Xu1993,TY1995,ET2020}
	and their references therein, where measurable selection theorems of multivalued 
	measurable functions were playing a key role in the study of random fixed points. 
	While Bharucha-Reid also posed random fixed point problems for random nonexpansive operators 
	in \cite[5.a, p653]{Bha1976}, the problems were deeply studied for a single random operator, 
	see, for example, \cite{Lin1988,Xu1990,Xu1993}. But 
	until  the present paper, we have not see even a common random fixed point theorem for a commutative 
	family of nonexpansive random operators since the usual method by measurable selection theorems fails 
	to work for this case. Let us explain the failure by taking the case of Theorem \ref{thm.1.18} 
	below for example. For each $T\in \mathcal{T}$, let $ F_{T}:\Omega\rightarrow 2^{V} $ defined by 
	$F_{T}(\omega)=\{v\in V:T(\omega,v)=v\} $, which is a nonempty closed set by the Kirk fixed point 
	theorem, so is $ \bigcap_{T\in \mathcal{T}}F_{T}(\omega) $ for each $ \omega\in \Omega $ by the  
	Lim's common fixed point theorem \cite{Lim1974}. When $V$ is separable, it is easy to see that 
	the graph $ Gr(F_{T}) $ of $ F_{T}$ $:=\{(\omega,x)\in \Omega\times B:T(\omega,x)-x=0\}$ is 
	$\mathcal{F}\otimes \mathcal{B}(V)$-measurable ($ \mathcal{B}(V) $ is the Borel $ \sigma $-algebra of 
	$ V $), but $ Gr(\bigcap_{T\in \mathcal{T}} F_{T}):=\bigcap_{T\in \mathcal{T}}Gr(F_{T})$ is not 
	necessarily $\mathcal{F}\otimes \mathcal{B}(V)$-measurable when $ \mathcal{T} $ is not countable. 
	So, the measurable selection theorems currently available in \cite{Wag1977} can not be used to obtain 
	a measurable selection of $\bigcap_{T\in \mathcal{T}}F_{T}$. Motivated by Guo et al. \cite{GZWY20}, 
	this paper provides a new method for solving the common random fixed point problem. Precisely speaking, we will lifting each random 
	nonexpansive operator $T:\Omega \times V\rightarrow V$ in $\mathcal{T}$ to a nonexpansive operator 
	$ \hat{T} $ from the closed $L^{0}$-convex subset $L^{0}(\mathcal{F},V)$ to itself, then 
	Theorem \ref{thm.1.18} will be an immediate corollary of Theorem \ref{thm.1.15}, 
	see the proof of Theorem \ref{thm.1.18} for details. Similarly, Theorems \ref{thm.1.19} and \ref{thm.1.20} are also respectively the immediate corollaries of Theorems \ref{thm.1.16} and \ref{thm.1.17}. The new method exhibits 
	the power of the idea of randomizing space theory. Following is a brief introduction 
	of the well known terminologies of random operators.

	\par Let $(\Omega,\mathcal{F},P)$ be a probability space and $(M,d)$ be a metric space.
	A mapping $X:\Omega\rightarrow M$ is said to be a random element \cite{Bha1976} (or, $\mathcal{F}$-random element) 
	if $X^{-1}(G):=\{\omega\in \Omega: X(\omega)\in G\}\in \mathcal{F}$ for each open set $G$ of $M$, furthermore, 
	a random element $X$ is said to be simple if it only takes finitely many values, and $X$ is said to be a strong 
	random element if $X$ is the pointwise limit of a sequence of simple random elements. It is well known that 
	$X:\Omega \rightarrow M$ is a strong random element if and only if both $X$ is a random element and its range $X(\Omega)$ is a separable subset of $M$.

	\par Let $(\Omega,\mathcal{F},P)$ be a probability space, $(M,d)$ and $(M_1,d_1)$ two metric spaces 
	and $T:\Omega\times M\rightarrow M_1$  a mapping. $T$ is called a random operator if $T(\cdot,x):\Omega\rightarrow M_1$ 
	is a random element for each $x\in M$. Further, $T$ is called a strong random operator if $T(\cdot,x)$ is a strong random 
	element for each $x\in M$. It is clear that each random operator $T:\Omega\times M\rightarrow M_1$ becomes a strong random 
	operator when $ M_{1} $ is separable.  $T$ is called a random nonexpansive operator (resp., random isometric operator) if 
	for each $\omega\in \Omega$, $d_{1}(T(\omega,x),T(\omega,y))\leq d(x,y)$ (resp., $d_{1}(T(\omega,x),T(\omega,y))= d(x,y)$) 
	for all $x,y\in M$. Two random operators $T_{1}:\Omega\times M\rightarrow M_1$ and $T_{2}:\Omega\times M\rightarrow M_1$ 
	are said to be commutative if for each $\omega\in \Omega$, $T_{1}(\omega,T_{2}(\omega,x))=T_{2}(\omega,T_{1}(\omega,x))$ for any $x\in M$.

	\par As an application of Theorem \ref{thm.1.15}, we can immediately obtain Theorem \ref{thm.1.18} below.

	\begin{thm}\label{thm.1.18}
		Let $( B , \| \cdot \|)$ be a Banach space over $\mathbb{K}$, $V$ be a weakly compact
		convex subset with normal structure of $B$ and  $\mathcal{T}$ a commutative family of
		strong random nonexpansive operators from $\Omega \times V$ to $V$.
		Then there is a strong random element $x^0(\cdot):\Omega
		\rightarrow V$ such that $T(\omega,x^0(\omega) )=x^0(\omega)$ for almost all
		$\omega \in \Omega$ and all $T\in \mathcal{T}$.
	\end{thm}

	\par Similarly, as a respective application of Theorems \ref{thm.1.16} and \ref{thm.1.17}, we can 
	immediately have Theorems \ref{thm.1.19} and \ref{thm.1.20} below.

	\begin{thm} \label{thm.1.19}
		Let $( B , \| \cdot \|)$ be a Banach space over $\mathbb{K}$, $V$ a weakly compact
		convex subset with normal structure of $B$ and $ T $ a strong random isometric operator from $\Omega \times V$ to $V$.
		Then there is a strong random element $x^0(\cdot):\Omega
		\rightarrow c(V,V)$ such that $T(\omega,x^0(\omega) )=x^0(\omega)$ for almost all
		$\omega \in \Omega$.
	\end{thm}

	\begin{thm} \label{thm.1.20}
			Let $( B , \| \cdot \|)$ be a Banach space over $\mathbb{K}$, $V$ a weakly compact
		convex subset with normal structure of $B$ and $\mathcal{T}$ a family of
		strong random surjective isometric operators from $\Omega \times V$ to $V$.
		Then there is a strong random element $x^0(\cdot):\Omega
		\rightarrow c(V,V)$ such that $T(\omega,x^0(\omega) )=x^0(\omega)$ for almost all
		$\omega \in \Omega$ and all $T\in \mathcal{T}$.
	\end{thm}

	\section{Proof of Theorem \ref{thm.1.14}}\label{sec.2}

	
	\par For the proof of Theorem \ref{thm.1.14}, let us first give Propositions \ref{prop.2.1} and \ref{prop.2.3} below, 
	which are concerned with the basic properties of random Chebyshev radius and  random Chebyshev center.

	\begin{prop} \label{prop.2.1}
		Let $E$, $G$ and $H$ be the same as in Definition \ref{defn. random Chebyshev center}. Then we have
		the following statements:
		\begin{enumerate}[(1)]
			\item $R(H,\cdot)$ is
			$\mathcal{T}_{\varepsilon,\lambda}$-continuous (namely, $R(H,\cdot):(E,
			\mathcal{T}_{\varepsilon,\lambda} )\rightarrow (L^0_{+}(\mathcal{F}),
			\mathcal{T}_{\varepsilon,\lambda})$ is continuous) and $L^0$-convex (namely, $R(H,\lambda x+(1-\lambda)y)\leq \lambda
			R(H,x)+(1-\lambda) R(H,y)$ for any $x,y\in E$ and any $\lambda \in
			L^0_{+}(\mathcal{F})$ with $0\leq \lambda \leq 1$).
			\item If $E$ is $\sigma$-stable, then $R(H,x)=R(\sigma(H),x )$ for any $x\in E$.
			\item $R(H,x)=R([Conv_{L^0}(H)]^{-}_{\varepsilon,\lambda},x)$ for any $x\in
			E$.
			\item  $\tilde{I}_{A} R(H,x)=R(\tilde{I}_{A}H,\tilde{I}_{A}x)$ for any $x\in
			E$ and any $A\in \mathcal{F}$. Further, if $E$ is $\sigma$-stable, then
			$$R( \sum_{n=1}^{\infty}\tilde{I}_{A_n}H_n,\sum_{n=1}^{\infty}\tilde{I}_{A_n}
			x_n)=\sum_{n=1}^{\infty}\tilde{I}_{A_n} R(H_n,x_n)$$ for any sequence $\{x_n,n\in
			\mathbb{N}\}$ in $E$, any sequence $\{H_{n},n\in \mathbb{N}\}$ of a.s. bounded
			nonempty subsets of $E$ and any $\{A_{n},n\in \mathbb{N}\} \in
			\Pi_{\mathcal{F}}$.
			\item If $E$ is $\sigma$-stable, then $R(H,G)=R(\sigma(H),\sigma(G) )$.
			\item $\tilde{I}_{A} R(H,G)=R(\tilde{I}_{A} H,\tilde{I}_{A} G)$ for any $A\in
			\mathcal{F}$. Further, if $E$ is $\sigma$-stable, then
			$$R(\sum_{n=1}^{\infty}\tilde{I}_{A_n}H_{n},
			\sum_{n=1}^{\infty}\tilde{I}_{A_n}G_{n})=\sum_{n=1}^{\infty}\tilde{I}_{A_n} R(H_{n},G_{n})$$ for any sequence $\{H_{n},n\in
			\mathbb{N}\}$ of a.s. bounded nonempty subsets of $E$,  any sequence $\{G_{n},n\in
			\mathbb{N}\}$ of  nonempty subsets of $E$  and any $\{A_{n},n\in
			\mathbb{N}\} \in \Pi_{\mathcal{F}}$.
			\item If $G$ is finitely stable and $C(H,G)\neq \emptyset$, then
			$\tilde{I}_{A} C(H,G)= C(\tilde{I}_{A}H,\tilde{I}_{A}G)$ for
			any $A\in \mathcal{F}$. Further, if $E$ is $\sigma$-stable, then 
			$$C(\sum_{n=1}^{\infty}\tilde{I}_{A_n}H_{n}, \sum_{n=1}^{\infty}\tilde{I}_{A_n}G_{n})=\sum_{n=1}^{\infty}\tilde{I}_{A_n}
			C( H_{n},G_{n})$$
			 for any sequence $\{H_{n},n\in \mathbb{N}\}$ of a.s.
			bounded nonempty subsets of $E$, any sequence $\{G_{n},n\in
			\mathbb{N}\}$ of  nonempty finitely stable subsets of $E$ and any $\{A_{n},n\in \mathbb{N}\} \in
			\Pi_{\mathcal{F}}$.
		\end{enumerate}
	\end{prop}
	\begin{proof}
		(1). It is clear, so is omitted.
		
		(2). It is obvious that $R(H,x)\leq R(\sigma(H),x)$. Conversely, for any
		sequence $\{x_n,n\in \mathbb{N}\}$ in $H$ and any $\{A_{n},n\in \mathbb{N}\} \in
		\Pi_{\mathcal{F}} $, $\|x-\sum_{n=1}^{\infty}\tilde{I}_{A_n}
		x_n\|=\sum_{n=1}^{\infty}\tilde{I}_{A_n} \|x- x_n\|\leq
		\sum_{n=1}^{\infty}\tilde{I}_{A_n} R(H,x)=R(H,x)$. Thus, $R(\sigma(H),x)\leq
		R(H,x)$.
		
		(3). First, we assert that $R(H,x)=R(Conv_{L^0}(H),x)$. It is obvious that
		$R(H,x)\leq R(Conv_{L^0}(H),x)$. On the other hand, for any $h'\in
		Conv_{L^0}(H)$, there exist $\{h_{i},i=1\sim n\}$ in $H$ and
		$\{\xi_{i},i=1\sim n\}$ in $L^0_{+}(\mathcal{F})$ with $\sum_{i=1}^{n}\xi_{i}=1$
		such that $h'=\sum_{i=1}^{n}\xi_{i} h_{i}$,
		then $\|x-h'\|=\|x-\sum_{i=1}^{n}\xi_{i} h_{i}\|\leq \sum_{i=1}^{n}\xi_{i} \|x-
		h_i\|\leq R(H,x)$, which implies that $R(Conv_{L^0}(H),x) \leq R(H,x)$. Then
		$R(Conv_{L^0}(H),x)= R(H,x)$.
		
		Second, we assert that $R(Conv_{L^0}(H),x)=
		R([Conv_{L^0}(H)]^{-}_{\varepsilon,\lambda},x)$. It is obvious that
		$R(Conv_{L^0}(H),x) \leq R([Conv_{L^0}(H)]^{-}_{\varepsilon,\lambda},x)$. On the
		other hand, for any $h\in [Conv_{L^0}(H)]^{-}_{\varepsilon,\lambda}$, there
		exists a sequence $\{h_n,n\in \mathbb{N}\}$ in $Conv_{L^0}(H)$ which converges
		to $h$. Then, we have
		\begin{align}\nonumber
			\|x-h\|&\leq \|x-h_n\|+\|h_n-h\|\nonumber\\
			&\leq R(Conv_{L^0}(H),x)+\|h_n-h\|\nonumber
		\end{align}
		for each $n\in \mathbb{N}$, which implies that
		$R([Conv_{L^0}(H)]^{-}_{\varepsilon,\lambda},x)\leq R(Conv_{L^0}(H),x)$. Thus,
		$R([Conv_{L^0}(H)]^{-}_{\varepsilon,\lambda},x)= R(Conv_{L^0}(H),x)$.
		
		To sum up, $R(H,x)=R([Conv_{L^0}(H)]^{-}_{\varepsilon,\lambda},x)$.
		
		(4). For any $x\in E$ and any $A\in \mathcal{F}$, we have
		\begin{align}\nonumber
			\tilde{I}_{A} R(H,x)&=\tilde{I}_{A}\bigvee\{\|x-h\|:h\in H\}\nonumber\\
			&=\bigvee\{\|\tilde{I}_{A} x-\tilde{I}_{A} h\|:h\in H\}\nonumber\\
			&=\bigvee\{\|\tilde{I}_{A} x- h'\|:h'\in \tilde{I}_{A}H\}\nonumber\\
			&=R(\tilde{I}_{A}H,\tilde{I}_{A}x).\nonumber
		\end{align}
		Further, since $\sum_{n=1}^{\infty}\tilde{I}_{A_n}
		x_n \in E$ and $\sum_{n=1}^{\infty}\tilde{I}_{A_n}H_{n}$ is still an a.s.
		bounded subset of $E$,  we have
		\begin{align}
			\tilde{I}_{A_n} R( \sum_{n=1}^{\infty}\tilde{I}_{A_n}
			H_{n},\sum_{n=1}^{\infty}\tilde{I}_{A_n}x_n )&= R(\tilde{I}_{A_n}
			H_{n},\tilde{I}_{A_n} x_n )\nonumber\\
			&=\tilde{I}_{A_n} R(H_{n},x_n)\nonumber	
		\end{align}
		for each $n\in \mathbb{N}$. Thus, $R(\sum_{n=1}^{\infty}\tilde{I}_{A_n}H_n ,
		\sum_{n=1}^{\infty}\tilde{I}_{A_n} x_n)=\sum_{n=1}^{\infty}\tilde{I}_{A_n}
		R(H_n,x_n)$.
		
		(5). It is obvious that $R(H,G)\geq R(H,\sigma(G))$. On the other hand, for
		any $x\in \sigma(G)$, there exist a sequence $\{x_n,n\in \mathbb{N}\}$ in $G$ and
		$\{A_{n},n\in \mathbb{N}\} \in \Pi_{\mathcal{F}}$ such that
		$x=\sum_{n=1}^{\infty}\tilde{I}_{A_n} x_n$, then
		\begin{align}\nonumber
			R(H,x)&=\sum_{n=1}^{\infty}\tilde{I}_{A_n} R(H,x_n) \geq R(H,G).\nonumber
		\end{align}
		It follows that $R(H,\sigma(G))\geq R(H,G)$, thus $R(H,G)= R(H,\sigma(G))$,
		which, combined with (2), implies that $R(H,G)=R(H,\sigma(G))=
		R(\sigma(H),\sigma(G))$.
		
		(6). For any $A\in \mathcal{F}$, we have
		\begin{align}\nonumber
			\tilde{I}_{A} R(H,G)&=\tilde{I}_{A} \bigwedge\{R(H,x): x\in G\}\nonumber\\
			&=\bigwedge\{ R(\tilde{I}_{A}H,\tilde{I}_{A}x): x\in G\}\nonumber\\
			&=\bigwedge\{ R(\tilde{I}_{A}H,x'): x'\in \tilde{I}_{A}G\}\nonumber\\
			&=R( \tilde{I}_{A}H,\tilde{I}_{A}G).\nonumber
		\end{align}
		Further, if $E$ is $\sigma$-stable,
		similar to the proof of (4), we have
		$$R(\sum_{n=1}^{\infty}\tilde{I}_{A_n}H_{n}, \sum_{n=1}^{\infty}\tilde{I}_{A_n}G_{n})=\sum_{n=1}^{\infty}\tilde{I}_{A_n} R(H_{n},G_{n}).$$
		
		(7). For any $A\in \mathcal{F}$ and any $x\in C(H,G)$,  we have
		$$R(\tilde{I}_{A}H,\tilde{I}_{A}x)=\tilde{I}_{A} R(H,x)=\tilde{I}_{A} R(H,G)=R(
		\tilde{I}_{A}H,\tilde{I}_{A}G),$$
		which implies that $\tilde{I}_{A}x\in C(
		\tilde{I}_{A}H,\tilde{I}_{A}G)$, and then $\tilde{I}_{A} C(H,G)\subset
		C(\tilde{I}_{A}H,\tilde{I}_{A}G)$. Conversely, for any $A\in
		\mathcal{F}$ and any $y\in C(\tilde{I}_{A}H,\tilde{I}_{A}G)$,
		$$\tilde{I}_{A} R(H,y)= R( \tilde{I}_{A}H,\tilde{I}_{A}y)=R(
		\tilde{I}_{A}H,\tilde{I}_{A}G)=\tilde{I}_{A} R(H,G). $$
		Take an arbitrary $x_0 \in C(H,G)$ and let
		$z=\tilde{I}_{A}y+\tilde{I}_{A^c}x_0$, then $z\in G$ and we have
		$R(H,z)=\tilde{I}_{A}R(H,y)+\tilde{I}_{A^c}R(H,x_0)=R(H,G)$, which implies that $z\in C(H,G)$. 
		Further, $y=\tilde{I}_{A}y=\tilde{I}_{A}z\in \tilde{I}_{A} C(H,G)$,
		then  $C(\tilde{I}_{A}H,\tilde{I}_{A}G) \subset  \tilde{I}_{A}
		C(H,G)$. Thus $\tilde{I}_{A}
		C(H,G)=C(\tilde{I}_{A}H,\tilde{I}_{A}G) $ for any $A\in
		\mathcal{F}$.
		
		By simple calculations as in (4), we can obtain 
		$$C(\sum_{n=1}^{\infty}\tilde{I}_{A_n}H_{n}, \sum_{n=1}^{\infty}\tilde{I}_{A_n}G_{n})=\sum_{n=1}^{\infty}\tilde{I}_{A_n}
		C( H_{n},G_{n}).$$
	\end{proof}

	\par The following Lemma \ref{lem.2.2}, established in \cite{GZWG20}, provides a useful property of a $\sigma$-stable set and is frequently used in this paper. Throughout this paper, $L^0_{++}(\mathcal{F})
	=\{\xi\in L^0(\mathcal{F}): \xi> 0 ~\text{on}~ \Omega\}$.

	\begin{lem}\label{lem.2.2}
		Let $G$ be a $\sigma$-stable subset of $L^0(\mathcal{F})$ such that $G$ has an
		upper (resp., lower) bound $\xi \in L^0(\mathcal{F})$, then for each
		$\varepsilon \in L^0_{++}(\mathcal{F})$ there exists some $g_{\varepsilon}\in G$
		such that $g_{\varepsilon}> \bigvee G-\varepsilon$ on $\Omega$ (resp.,
		$g_{\varepsilon}< \bigwedge G+\varepsilon$ on $\Omega$).
	\end{lem}

	\par As shown by \cite[Lemma 2.5]{GWXYC24}, if $G$ is a finitely stable and $\mathcal{T}_{\varepsilon,\lambda}$-closed 
	subset of some $\sigma$-stable subset of an $RN$ module $(E,\|\cdot\|)$, then $G$ is $\sigma$-stable. Clearly, 
	if $(E,\|\cdot\|)$ is  $\mathcal{T}_{\varepsilon,\lambda}$-complete, then  $E$ is $\sigma$-stable.

	\begin{prop} \label{prop.2.3}
		Let $( E , \| \cdot \|)$ be a $\mathcal{T}_{\varepsilon,\lambda}$-complete
		$RN$ module over $\mathbb{K}$ with base $(\Omega,\mathcal{F},P)$, $G$
		an $L^0$-convexly compact subset of $E$ and $H$ an a.s. bounded nonempty subset
		of $E$. Then $C(H,G)$ is a nonempty
		$\mathcal{T}_{\varepsilon,\lambda}$-closed $L^0$-convex subset of $G$.
	\end{prop}
	\begin{proof}
		For each $n\in \mathbb{N}$, let $C_n=\{x\in G : R(H,x) \leq R(H,G)+ \frac{1}{n}\}$.
		First, since $G$ is $\sigma$-stable, then $\{R(H,x): x\in G\}$ is a
		$\sigma$-stable subset of $L^0_{+}(\mathcal{F})$ by (4) of Proposition
		\ref{prop.2.1}, and hence $C_n$ is nonempty by Lemma \ref{lem.2.2}. Second, since $G$
		is $\mathcal{T}_{\varepsilon,\lambda}$-closed and $R(H,\cdot)$ is
		$\mathcal{T}_{\varepsilon,\lambda}$-continuous, then $C_n$ is
		$\mathcal{T}_{\varepsilon,\lambda}$-closed. Finally, for any $x_1,x_2\in C_n$
		and any $\lambda \in L^0_{+}(\mathcal{F})$ with $0 \leq\lambda \leq 1$, 
		$\lambda x_1+(1-\lambda)x_2 \in G$ and we have
		\begin{align}\nonumber
			R(H,\lambda x_1+(1-\lambda)x_2)&\leq \lambda
			R(H,x_1)+(1-\lambda)R(H,x_2)\nonumber\\
			&\leq R(H,G)+ \frac{1}{n},\nonumber
		\end{align}
		thus $C_n$ is $L^0$-convex. Therefore, $\{C_n,n\in \mathbb{N}\}$ is a family of nonempty
		$\mathcal{T}_{\varepsilon,\lambda}$-closed $L^0$-convex subset of $G$. Further, it is
		obvious that $\{C_n,n\in \mathbb{N}\}$ has finite intersection property. By the
		$L^0$-convex compactness of $G$, $C(H,G)=\bigcap_{n=1}^{\infty}
		C_n\neq \emptyset $.
	\end{proof}

	\par It is straightforward to check that complete normal structure implies normal structure for any  bounded closed convex subset of a  Banach space by taking $W_{\alpha}=W$ for any $\alpha\in \Lambda$ in Definition \ref{defn.cns}. However, as shown in Proposition \ref{prop.2.4} below, this is not so straightforward in the case of an $RN$ module.

	\begin{prop}\label{prop.2.4}
		Let $( E , \| \cdot \|)$ be a $\mathcal{T}_{\varepsilon,\lambda}$-complete
		$RN$ module over $\mathbb{K}$ with base $(\Omega,\mathcal{F},P)$ and
		$G$ be an a.s. bounded $\mathcal{T}_{\varepsilon,\lambda}$-closed $L^0$-convex
		subset of $E$ such that $G$ has random complete normal structure. Then $G$ has
		random normal structure.
	\end{prop}
	\begin{proof}
		Assume that $G$ does not have random normal structure, then there exists a
		nonempty $\mathcal{T}_{\varepsilon,\lambda}$-closed $L^0$-convex subset $W$ of
		$G$ , with $W$ containing more than one point,  such that for any $x\in W$, there exists
		$A_{x}\subset (D(W)>0)$ with $P(A_{x})>0$ such that
		\begin{equation}\label{1234}
			\tilde{I}_{A_x}R(W,x)=\tilde{I}_{A_x} D(W).
		\end{equation}
		
		Let $\Lambda$ be a $B_{\mathcal{F}}$-stable directed set and $W_{\alpha}=W$
		for any $\alpha \in \Lambda$, then $\{W_\alpha, \alpha \in \Lambda\}$ is a
		decreasing consistent net of $\sigma$-stable subsets of $W$ and $R(W_\alpha,W)$
		$=R(W,W)$ for any $\alpha \in \Lambda$. Since $G$ has random complete normal
		structure, then $[C(W,W) ]^{-}_{\varepsilon,\lambda}\neq \emptyset$.
		For any given $x_0\in C(W,W)$, by (\ref{1234}),
		there exists $A_{x_0}\subset (D(W)>0)$ with $P(A_{x_0})>0$ such that
		\begin{equation}\label{12345}
			\tilde{I}_{A_{x_0}}R(W,x_0)=\tilde{I}_{A_{x_0}}D(W).
		\end{equation}
		For $A_{x_0}$, there exists $y_0\in W$ such that $\tilde{I}_{A_{x_0}} y_0
		\notin \tilde{I}_{A_{x_0}} [C(W,W)]^{-}_{\varepsilon,\lambda}$. It
		follows that $\tilde{I}_{A_{x_0}}y_0 \notin
		\tilde{I}_{A_{x_0}}C(W,W)=C(\tilde{I}_{A_{x_0}}W,
		\tilde{I}_{A_{x_0}}W)$, which implies that $R(\tilde{I}_{A_{x_0}}W,\tilde{I}_{A_{x_0}}y_0)\neq R(\tilde{I}_{A_{x_0}}W,
		\tilde{I}_{A_{x_0}}W)$.
		Thus, there exists $C\subset A_{x_0}$ with $P(C)>0$ such
		that
		\begin{align}
			\tilde{I}_{C} R(\tilde{I}_{A_{x_0}}W,\tilde{I}_{A_{x_0}}y_0) &>\tilde{I}_{C} R(\tilde{I}_{A_{x_0}}W,\tilde{I}_{A_{x_0}}W)\nonumber\\
			&=\tilde{I}_{C}R(W,W)\nonumber\\
			&=\tilde{I}_{C}R(W,x_0)\nonumber\\
			&=\tilde{I}_{C}D(W)\nonumber
		\end{align}
		on $C$, which  contradicts  with $\tilde{I}_{C} R(\tilde{I}_{A_{x_0}}W,\tilde{I}_{A_{x_0}}y_0)=\tilde{I}_{C} R(W,y_0)\leq \tilde{I}_{C}D(W)$.
	\end{proof}
	
  \par To proceed,	we need the notion of random asymptotic centers as follows.
	
	\begin{defn}\label{defn.2.5}
		Let $( E , \| \cdot \|)$ be an
		$RN$ module over $\mathbb{K}$ with base $(\Omega,\mathcal{F},P)$, $G$
		a nonempty subset of $E$ and  $\mathcal{E}=\{B_\alpha, \alpha\in
		\Lambda\}$ a decreasing net of a.s. bounded nonempty subsets of $E$. Define $AR(\mathcal{E},\cdot):E\rightarrow L^{0}_{+}(\mathcal{F})$ by
		$$AR(\mathcal{E},x)=\bigwedge \{R(B_\alpha,x): \alpha\in \Lambda \},~\forall x\in E.$$
		Then
		$AR(\mathcal{E},G)=\bigwedge\{AR(\mathcal{E},x): x\in G\} $
		and
		$AC(\mathcal{E},G)=\{x\in G:
		AR(\mathcal{E},x)=AR(\mathcal{E},G)  \} $
		are called the
		random asymptotic radius and random asymptotic center of $\mathcal{E}$ respect
		to $G$, respectively.
	\end{defn}

	\par For an $RN$ module $(E,\|\cdot\|)$ over $\mathbb{K}$ with base $(\Omega,\mathcal{F},P)$, 
	$d: E\times E\rightarrow L^0_{+}(\mathcal{F})$ defined by
	$d(x,y)=\|x-y\|$ for any $x,y\in E$, is clearly a random metric on $E$ and thus
	$(E,d)$ is a random metric space \cite{Guo92}. Let $CB(E)$ be the family of nonempty a.s.
	bounded and $\mathcal{T}_{\varepsilon,\lambda}$-closed subsets of $E$ and
	$CB_{\sigma}(E)=\{G \in CB(E): G ~\text{is}~ \sigma\text{-stable}  \}$. Define
	the random Hausdorff metric
	$H: CB_{\sigma}(E) \times CB_{\sigma}(E) \rightarrow L^{0}_{+}(\mathcal{F})$ by
	$$H(G_1,G_2)=\max\{ \bigvee_{x\in G_1}d(x, G_2),  \bigvee_{x_2\in G_2}d(x_2,
	G_1)  \}$$ for any $G_1$ and $G_1$ in $CB_{\sigma}(E)$, where
	$d(x,G)=\bigwedge\{d(x,g):g\in G\}$ denotes the random distance from $x\in E$ to a
	nonempty subset $G$ of $E$. Then, $(CB_{\sigma}(E), H)$ is a random metric space and
	the $L^0$-topology on $CB_{\sigma}(E)$ is  denoted by $\mathcal{T}^{H}_{c}$. Here, we do not give the notion of
	the $L^0$-topology,  see \cite{GMT24,GWYZ20} for the detailed definition of
	the $L^0$-topology, it suffices to say that a net $ \{G_{\alpha},\alpha\in \Lambda\} $ in $CB_{\sigma}(E)$ converges in $\mathcal{T}^{H}_{c}$
	to $ G $ iff for any $ \varepsilon\in L^{0}_{++}(\mathcal{F}) $ there exists $ \alpha_{\varepsilon}\in \Lambda $ such that $ H(G_{\alpha},G)<\varepsilon $ on $ \Omega $ whenever $ \alpha\geq \alpha_{\varepsilon} $.

	\par We would like to remind readers that an $L^0$-convexly compact subset $G$ of an $RN$ module 
	must be a.s. bounded (see \cite[Lemma 2.19]{GZWW20} or \cite[Proposition 3.3]{WZZ22} for details).

	\par Lemmas \ref{lem.2.6} and \ref{lem.2.7} below are crucial for the proof of Theorem
	\ref{thm.1.14}.
	\begin{lem}\label{lem.2.6}
		Let $( E , \| \cdot \|)$ be a $\mathcal{T}_{\varepsilon,\lambda}$-complete
		$RN$ module over $\mathbb{K}$ with base $(\Omega,\mathcal{F},P)$, $W$
		an $L^0$-convexly compact subset of $E$ and $\mathcal{E}=\{W_\alpha, \alpha \in
		\Lambda\}$ a decreasing consistent net of $\sigma$-stable subsets of $W$ with
		 $R(W_\alpha,W)=R(W,W)$ for each $\alpha \in \Lambda$.
		Then we have the following statements:
		\begin{enumerate}[(1)]
			\item $AR(\mathcal{E},x)=R(W,W)$ for any $x\in [\bigcup_{\alpha \in
				\Lambda} C(W_\alpha,W)]^{-}_{\varepsilon,\lambda}$.
			\item If $AR(\mathcal{E},x)=0$, then
			$\{[W_{\alpha}]^{-}_{\varepsilon,\lambda}, \alpha \in \Lambda \}$ converges in
			$\mathcal{T}^{H}_{c}$ to $\{x\}$.
		\end{enumerate}
	\end{lem}
	
	\begin{proof}
		(1). For any $x\in \bigcup_{\alpha \in \Lambda} C(W_\alpha,W)$,
		there exists $\beta \in \Lambda$ such that $x\in  C(W_{\beta},W)$,
		then we have
		\begin{equation}\nonumber
			AR(\mathcal{E},x)\leq R(W_{\beta},x)=R(W_{\beta},W)=R(W,W).
		\end{equation}
		On the other hand, for any $x\in \bigcup_{\alpha \in \Lambda} C(W_\alpha,W)$, we have
		\begin{equation}\nonumber
			R(W,W)=R(W_{\alpha},W)\leq R(W_{\alpha},x) ~\forall~ \alpha \in \Lambda,
		\end{equation}
		which implies that $R(W,W) \leq AR(\mathcal{E},x)$. Therefore,
		$AR(\mathcal{E},x)=R(W,W)$ for any $x \in \bigcup_{\alpha \in \Lambda}
		C(W_\alpha,W)$. By the $ \mathcal{T}_{\varepsilon,\lambda}$-continuity of $AR(\mathcal{E},\cdot)$ (see (1) of Lemma \ref{lem.4.3} below), we
		complete the proof.
		
		(2). Since $\Lambda$ is a $B_{\mathcal{F}}$-stable directed set, it is easy to check that $\{R(W_{\alpha},x): \alpha \in \Lambda\}$ is a
		$\sigma$-stable subset of $L^0_{+}(\mathcal{F})$ by (4) of Proposition
		\ref{prop.2.1}. By Lemma \ref{lem.2.2}, for any
		$\varepsilon \in L^0_{++}(\mathcal{F})$ there exists $\alpha_{0}\in \Lambda$
		such that
		$R(W_{\alpha},x)<AR(\mathcal{E},x)+\varepsilon=\varepsilon$ on $\Omega$ for any $\alpha \in \Lambda$
		with $\alpha \geq \alpha_0$. Further, for any $\alpha\in \Lambda$,
		\begin{align}\nonumber
			H(\{x\},[W_{\alpha}]^{-}_{\varepsilon,\lambda})
			&=\max\{ \bigvee_{x\in \{x\}}d(x, [W_{\alpha}]^{-}_{\varepsilon,\lambda}),
			\bigvee_{y\in [W_{\alpha}]^{-}_{\varepsilon,\lambda}}d(y, \{x\})  \} \\
			\nonumber
			&= R([W_{\alpha}]^{-}_{\varepsilon,\lambda},x)  \\ \nonumber
			&= R(W_{\alpha},x).  \nonumber
		\end{align}
		Then $H(\{x\},[W_{\alpha}]^{-}_{\varepsilon,\lambda})<\varepsilon$ on $\Omega$
		for any $\alpha \in \Lambda$ with $\alpha \geq \alpha_0$, which implies that
		$\{[W_{\alpha}]^{-}_{\varepsilon,\lambda}, \alpha \in \Lambda \}$ converges in
		$\mathcal{T}^{H}_{c}$ to $\{x\}$.
	\end{proof}

	\begin{lem}\label{lem.2.7}
		Let $( E , \| \cdot \|)$ be a $\mathcal{T}_{\varepsilon,\lambda}$-complete
		$RN$ module over $\mathbb{K}$ with base $(\Omega,\mathcal{F},P)$ and $G$ be
		a nonempty $\mathcal{T}_{\varepsilon,\lambda}$-closed $L^0$-convex subset of $E$. If there exists  
		a sequence $\{x_n,n\in \mathbb{N}\}$  in $G$ such that for some $c\in
		L^{0}_{++}(\mathcal{F})$, $\|x_n-x_m\|\leq c$, $\|x_{n+1}-\overline{x}_n\|\geq
		c-\frac{1}{n^2}$ for all $n\geq 1$, $m\geq 1$, where
		$\overline{x}_n=\frac{1}{n}\sum_{i=1}^{n}x_i$, then $G$ does not have random
		normal structure.
	\end{lem}
	
	\begin{proof}
		Let $H=[ Conv_{L^0}(\{x_n,n\in \mathbb{N}\}) ]^{-}_{\varepsilon,\lambda}$. For
		any fixed $\hat{x}\in Conv_{L^0}(\{x_n,n\in \mathbb{N}\})$, let 
		$\hat{x}=a_1x_1+a_2x_2+\cdots+ a_lx_l$ with $a_i\in L^{0}_{+}(\mathcal{F})$ for any
		$i=1,\cdots,l$ and $\sum_{i=1}^{l}a_i=1$, $a:=\max\{a_1,
		a_2,\cdots, a_l\}$ and $a_{l+1}=\cdots=a_{n}=0$ for any $n>l$, we have
		\begin{align}\nonumber
			\|x_{n+1}-\hat{x}\|
			&=\| nax_{n+1}-a \sum_{i=1}^{n}x_{i} -nax_{n+1}+a \sum_{i=1}^{n}x_{i} +
			\sum_{i=1}^{n}a_{i} x_{n+1} - \sum_{i=1}^{n}a_{i} x_{i}    \| \nonumber  \\
			&=\| na(x_{n+1}-\frac{1}{n} \sum_{i=1}^{n}x_{i}) - \sum_{i=1}^{n}
			(a-a_i)(x_{n+1}-x_i)   \| \nonumber \\
			&\geq na \| x_{n+1}-\overline{x}_{n} \| -  \sum_{i=1}^{n} (a-a_i) \|
			x_{n+1}-x_i \| \nonumber \\
			&\geq na (c-\frac{1}{n^2})- (na-1)c \nonumber \\
			& \geq c-\frac{1}{n}, \nonumber
		\end{align}
		which implies that
		\begin{align}
			R(\{x_n,n\in \mathbb{N}\},\hat{x} )&\geq
			R(\{x_{n+1},n>l\},\hat{x} )\nonumber\\
			&=\bigvee \{\|\hat{x}-x_{n+1}\|:n>l\}\nonumber\\
			&\geq\bigvee\{c-\frac{1}{n},n>l\}\nonumber\\
			&=c.\nonumber	
		\end{align}
		On the other hand, it is easy to check that $R(
		\{x_n,n\in \mathbb{N}\},\hat{x} )\leq c$.
		Thus, $R( \{x_n,n\in \mathbb{N}\},\hat{x} )=c$.
		Further, (3) of Proposition
		\ref{prop.2.1} shows that $R( H,\hat{x} )=c$. By the
		$\mathcal{T}_{\varepsilon,\lambda}$-continuity of $R(H,\cdot)$, we have
		$R(H,x)=c$ for any $x\in H$,  and hence $D(H)=c$. Then, the a.s bounded
		$\mathcal{T}_{\varepsilon,\lambda}$-closed $L^0$-convex subset $H$ of $G$ does
		not have a nondiametral point, namely, $G$ does not have random normal
		structure.
	\end{proof}

	\par For $\xi \in L^0( \mathcal{F} ,\mathbb{K} )$ with a representative $\xi^0$, $(\xi^0
	)^{-1}:\Omega\rightarrow \mathbb{K}$ is defined by $(\xi^0 )^{-1}(\omega)=\frac{1}{\xi^0(\omega)}$ if $\xi^0(\omega)\neq
	0$ and $0$ otherwise. Then the equivalence class
	of $(\xi^0)^{-1}$ is called the generalized inverse of $\xi$, denoted by $\xi^{-1}$. It is
	clear that $\xi \cdot \xi^{-1}=\tilde{I}_{(\xi\neq 0)}$.
	
	\par Now, we are ready to prove Theorem \ref{thm.1.14}.

	\begin{proof}[\textbf{Proof of Theorem \ref{thm.1.14}}]
		$\mathbf{Necessity}$. It is obvious by Proposition \ref{prop.2.4}.
		
		$\mathbf{Sufficiency}$. Suppose that $G$ does not have random complete normal
		structure, then there exists a nonempty
		$\mathcal{T}_{\varepsilon,\lambda}$-closed $L^0$-convex subset $W$ of $G$, with $W$
		containing more than one point,  together with a decreasing consistent net $\mathcal{E}:=\{W_\alpha, \alpha \in \Lambda\}$
		of $\sigma$-stable subsets of $W$
		satisfying $R(W_\alpha,W)=R(W,W)$ for each $\alpha \in \Lambda$, but there
		exists some $B\subset (D(W)>0)$ with $P(B)>0$ such that $\tilde{I}_{B}y \in
		\tilde{I}_{B} [\bigcup_{\alpha \in \Lambda} C(W_\alpha,W)]
		^{-}_{\varepsilon,\lambda} $ for any $y\in W$ (note that $[\bigcup_{\alpha \in
			\Lambda} C( W_\alpha,W)]^{-} _{\varepsilon,\lambda} \neq \emptyset$ by
		Proposition \ref{prop.2.3}). Consequently,
		$$\tilde{I}_{B} [\bigcup_{\alpha \in
			\Lambda} C(W_\alpha,W)]^{-}_{\varepsilon,\lambda}= \tilde{I}_{B} W.$$
		Without loss of generality, we can assume that $B=\Omega$ (otherwise, we
		consider the $RN$ module $(E_{B},\|\cdot\|_{E_{B}})$ over $\mathbb{K}$ with base
		$(B,\mathcal{F}_{B}, P_{B})$, where $E_{B}=\tilde{I}_{B}E$, $\|\cdot\|_{E_{B}}$
		be the restriction of $\|\cdot\|$ to $E_{B}$,  $\mathcal{F}_{B}=\{A\cap B:A\in
		\mathcal{F}\}$ and $P_{B}(A\cap B)=\frac{P(A\cap B)}{P(B)}$ for any $A\cap B\in
		\mathcal{F}_{B}$). Then $D(W)\in L^{0}_{++}(\mathcal{F})$ and $$ [\bigcup_{\alpha
			\in \Lambda} C(W_\alpha,W)]^{-}_{\varepsilon,\lambda}=W.$$
		
		First, we prove that $R(W,W)>0$ on $\Omega$. Otherwise, there exists $C\in
		\mathcal{F}$ with $P(C)>0$ such that $R(W,W)=0$ on $C$.
		Since $D(W)\in L^{0}_{++}(\mathcal{F})$ and  $\{\|x-y\|:x,y\in W\}$ is a $\sigma$-stable subset of
		$L^{0}_{+}(\mathcal{F})$, by Lemma \ref{lem.2.2}, for $\varepsilon=\frac{D(W)}{2}\in
		L^{0}_{++}(\mathcal{F})$, there exist $x_0, y_0\in W$ such that
		\begin{align}\nonumber
			\|x_0-y_0\|> D(W)-\varepsilon=\frac{D(W)}{2} ~\text{on}~ \Omega.
		\end{align}
		Further, since $ [\bigcup_{\alpha
			\in \Lambda} C(W_\alpha,W)]^{-}_{\varepsilon,\lambda}=W$, by (1) of Lemma \ref{lem.2.6}, we have
		$AR(\mathcal{E},x_0)=AR(\mathcal{E},y_0)=R(W,W)=0$ on $C$. It follows that
		$\{ \tilde{I}_{C} [W_{\alpha}]^{-}_{\varepsilon,\lambda}, \alpha \in \Lambda \}$
		converges in $\mathcal{T}^{H}_{c}$ to two distinct sets $\{\tilde{I}_{C}x_0\}$
		and $\{\tilde{I}_{C}y_0\}$ by (2) of Lemma \ref{lem.2.6}, which contradicts with
		the fact that $\mathcal{T}^{H}_{c}$ is a Hausdorff topology.
		
		\par Second, we prove that $G$ does not have random normal structure, which leads to a 
		contradiction with the assumption of sufficiency. By Lemma
		\ref{lem.2.7}, it suffices to show that there exists a sequence $\{x_{n},n\in
		\mathbb{N}\}$ in $W$ such that for some $c\in L^{0}_{++}(\mathcal{F})$,
		$\|x_n-x_m\|\leq c$, $\|x_{n+1}-\overline{x}_n\|\geq c-\frac{1}{n^2}$ for all
		$n\geq 1$, $m\geq 1$, where $\overline{x}_n=\frac{1}{n}\sum_{i=1}^{n}x_i$.
		
		We construct $\{x_{n}, n \in \mathbb{N}\}$ by using the method of induction.
		Let $x_1$ be an arbitrary element of $\bigcup_{\alpha \in \Lambda} C(
		W_\alpha,W)$, there exists $\alpha_{1}\in \Lambda$ such that
		$x_1\in C( W_{\alpha_{1}},W)$, then $R(W_{\alpha_{1}},x_1)=
		R(W_{\alpha_{1}},W)=R(W,W)$. By Lemma \ref{lem.2.2}, there exists $x_2\in
		W_{\alpha_{1}}$ such that $\|x_1-x_2\|\geq R(W_{\alpha_{1}},x_1)-1= R(W,W)-1$
		and $\|x_1-x_2\|\leq R(W_{\alpha_{1}},x_1) = R(W,W)$.
		
		\par Suppose now $x_1,x_2,\cdots,x_n$ have been chosen in $W$
		such that  $\|x_i-x_j\|\leq R(W,W)$ and $\|x_{i+1}-\overline{x}_{i}\|\geq
		R(W,W)-\frac{1}{i^2}$ for all $1\leq i\leq n-1$ and $1\leq j\leq n$, where
		$\overline{x}_{i}=\frac{1}{i}\sum_{k=1}^{i}x_k$. We proceed to choose $x_{n+1}$
		as follows. Since $\mathcal{E}=\{W_\alpha, \alpha \in \Lambda\}$ is a decreasing consistent net and
		$R(W,W)=AR(\mathcal{E},x_i)=AR(\mathcal{E},\overline{x}_{n})$ for any
		$i\in\{1,2,\cdots,n\}$ by (1) of Lemma \ref{lem.2.6}. By Lemma \ref{lem.2.2}, there exists $\alpha_{n} \in \Lambda$ such that
		\begin{align}\label{1.2}
			R(W_{\alpha_{n}},x_i) &\leq AR(\mathcal{E},x_i)+\frac{1}{n^2 (n+1)}\nonumber\\
			&=R(W,W)+\frac{1}{n^2 (n+1)}
		\end{align}
		for any $i\in\{1,2,\cdots,n\}$ and
		\begin{align}
			R(W_{\alpha_{n}},\overline{x}_{n}) &\leq AR(\mathcal{E},\overline{x}_{n})+\frac{1}{n^2 (n+1)}\nonumber\\
			&=R(W,W)+\frac{1}{n^2 (n+1)}.\nonumber	
		\end{align}
		Since $\{\|\overline{x}_{n}-y\|:y\in W_{\alpha_{n}}\}$ is a $\sigma$-stable subset
		of $L^0_{+}(\mathcal{F})$, by Lemma \ref{lem.2.2}, there exists $z_0\in W_{\alpha_{n}}$ such
		that
		\begin{align}\label{1.3}
			\|\overline{x}_{n}-z_0\|&\geq R(W_{\alpha_{n}},\overline{x}_{n})- \frac{1}{n^2
				(n+1)}\nonumber \\
			&\geq AR(\mathcal{E},\overline{x}_{n})- \frac{1}{n^2 (n+1)} \nonumber \\
			&= R(W,W)- \frac{1}{n^2 (n+1)}.
		\end{align}
		Let $z_i=t_i x_i+ (1-t_i)z_{i-1}$, where
		$$t_i=\max \{ \|x_i-z_{i-1}\|^{-1}
		(\|x_i-z_{i-1}\|-R(W,W)), 0\},$$
		$i=1,2,\cdots,n$.
		Let $x_{n+1}:=z_n$, we assert that $\|x_{n+1}-x_i\|\leq R(W,W)$ for any
		$i=1,2,\cdots,n$ and $\|x_{n+1}-\overline{x}_n\|\geq R(W,W)-\frac{1}{n^2}$.
		
		Indeed, first, for any $i\in \{1,\cdots,n\}$,
		\begin{align}
			\|z_{i}-x_{i}\|&=\|t_i x_i+ (1-t_i)z_{i-1}-x_{i}\|\nonumber\\
			&=(1-t_{i})\|z_{i-1}-x_{i}\|\nonumber\\
			&=(1-\max \{ \|x_i-z_{i-1}\|^{-1}
			(\|x_i-z_{i-1}\|-R(W,W)), 0\})\|z_{i-1}-x_{i}\|\nonumber\\
			&=\min\{R(W,W),\|z_{i-1}-x_{i}\|\}\nonumber\\
			&\leq R(W,W),\nonumber
		\end{align}
		then we have
		\begin{align}
			\|x_{n+1}-x_i\|&=\|z_{n}-x_{i}\|\nonumber\\
			&=\|t_{n}x_{n}+(1-t_{n})z_{n-1}-x_{i}\|\nonumber\\
			&\leq t_{n}\|x_{n}-x_{i}\|+(1-t_{n})\|z_{n-1}-x_{i}\|\nonumber\\
			&\quad\vdots\nonumber\\
			&\leq t_{n}\|x_{n}-x_{i}\|+(1-t_{n})t_{n-1}\|x_{n-1}-x_{i}\|+\cdots\nonumber\\
			&\quad\quad\quad+(1-t_{n})(1-t_{n-1})\cdots(1-t_{i+2})t_{i+1}\|x_{i+1}-x_{i}\|\nonumber\\
			&\quad\quad\quad+(1-t_{n})(1-t_{n-1})\cdots(1-t_{i+2})(1-t_{i+1})\|z_{i}-x_{i}\|\nonumber\\
			&\leq R(W,W).\nonumber
		\end{align}
		
		Second, (\ref{1.2}) implies that
		$$\|z_0-x_i\| \leq  R(W,W)+\frac{1}{n^2 (n+1)},$$
		then we have
		\begin{align}
			\|z_{i-1}-x_{i}\|&=\|t_{i-1}x_{i-1}+(1-t_{i-1})z_{i-2}-x_{i}\|\nonumber\\
			&\leq t_{i-1}\|x_{i-1}-x_{i}\|+(1-t_{i-1})\|z_{i-2}-x_{i}\|\nonumber\\
			&\quad \vdots\nonumber\\
			&\leq t_{i-1}\|x_{i-1}-x_{i}\|+(1-t_{i-1})t_{i-2}\|x_{i-2}-x_{i}\|+\cdots\nonumber\\
			&\quad\quad\quad+(1-t_{i-1})(1-t_{i-2})\cdots(1-t_{2})t_{1}\|x_{1}-x_{i}\|\nonumber\\
			&\quad\quad\quad+(1-t_{i-1})(1-t_{i-2})\cdots(1-t_{1})\|z_{0}-x_{i}\|\nonumber\\
			&\leq R(W,W)+\frac{1}{n^2 (n+1)}.\nonumber	
		\end{align}
		It follows that
		\begin{align}
			\|x_{n+1}-z_0\|&=\|z_n-z_0\| \nonumber \\
			&\leq \sum _{i=1}^{n} \|z_i-z_{i-1}\| \nonumber \\
			&= \sum _{i=1}^{n} t_i \|x_i-z_{i-1}\|  \nonumber \\
			&=\sum _{i=1}^{n} \max \{ \|x_i-z_{i-1}\|-R(W,W)  , 0\}\nonumber \\
			&\leq \frac{n}{n^2 (n+1)},\nonumber
		\end{align}
		which, combined with (\ref{1.3}), implies that
		\begin{align}\nonumber
			\|x_{n+1}-\overline{x}_n\| &\geq \|\overline{x}_n -z_0\| - \|z_0- x_{n+1}\|
			\geq R(W,W)- \frac{1}{n^2 }. \nonumber
		\end{align}
			
	\end{proof}
	

	\par Let $(E,\|\cdot\|)$ be an $RN$ module over $\mathbb{K}$ with base
	$(\Omega,\mathcal{F},P)$ and $\xi=\bigvee\{\|x\|:x\in E\}$ (generally,
	$\xi\in \bar{L}^{0}_{+}(\mathcal{F}):=\{\xi\in \bar{L}^0(\mathcal{F}): \xi\geq 0\}$). $(E,\|\cdot\|)$ is said to have full
	support if $P((\xi>0))=1$.
	
	\par In the remainder of this paper, all $RN$ modules mentioned are assumed to have full support. Further, we employ the
	following notations for a brief introduction to random uniformly convex $RN$ modules:\\
	{\small $\varepsilon_{\mathcal{F}}[0,2]=\{ \varepsilon \in
		L^{0}_{++}(\mathcal{F})| $ there exists a positive number$ ~\lambda~ $such that$
		~\lambda\leq \varepsilon \leq 2\}$.}
	{\small $\delta_{\mathcal{F}}[0,1]=\{ \delta \in L^{0}_{++}(\mathcal{F})| $
		there exists a positive number$ ~\eta~ $such that$ ~\eta \leq \delta \leq 1\}$.}
	$A_{x}=(\|x\|>0)$, $A_{xy}=A_{x}\cap A_{y}$ and $B_{xy}=A_{x}\cap A_{y} \cap
	A_{x-y}$ for any $x$ and $y$ in $E$.
	
	\par An $RN$ module $( E , \| \cdot \|)$ is said to be random uniformly convex \cite{GZ10} if
	for each $\varepsilon \in \varepsilon_{\mathcal{F}}[0,2]$, there exists $\delta
	\in \delta_{\mathcal{F}}[0,1]$ such that $\|x-y\|\geq \varepsilon$ on $D$ always
	implies $\|x+y\|\leq 2(1-\delta)$ on $D$ for any $x,y \in U(1)$ and any $D\in
	\mathcal{F}$ such that $D\subset B_{xy}$ with $P(D)>0$, where $U(1)=\{z\in E:
	\|z\| \leq 1\}$, called the random closed unit ball of $E$.

	\begin{cor}\label{cor.2.8}
		Let $( E , \| \cdot \|)$ be a $\mathcal{T}_{\varepsilon,\lambda}$-complete
		random uniformly convex $RN$ module over $\mathbb{K}$ with base
		$(\Omega,\mathcal{F},P)$ and $G$ be a nonempty a.s. bounded
		$\mathcal{T}_{\varepsilon,\lambda}$-closed and $L^0$-convex subset of $E$. Then
		$G$ has random complete normal structure.
	\end{cor}
	
	\begin{proof}
		By \cite[Theorem 4.6]{GZ10}, $E$ is random reflexive, it follows from \cite[Corollary 2.23]{GZWW20} that $G$ is $L^0$-convexly 
		compact. Furthermore, \cite[Theorem 3.6]{GZWG20} shows 
		that $G$ has random normal structure. Therefore, $G$ has random complete normal structure by Theorem \ref{thm.1.14}.
	\end{proof}

	\par Let $B$ be a Banach space over $\mathbb{K}$, $ V $ a closed convex subset of $ B $ and $L^0(\mathcal{F},B)$ (resp., $L^0(\mathcal{F},V)$) the set of equivalence
	classes of strong random elements from $(\Omega,\mathcal{F},P)$ to
	$B$ (resp., $V$). It is known from \cite{Guo10} that $L^0(\mathcal{F},B)$ (resp., $L^0(\mathcal{F},V)$) is a $\mathcal{T}_{\varepsilon,\lambda}$-complete $RN$ module (resp., $\mathcal{T}_{\varepsilon,\lambda}$-closed $L^{0}$-convex subset of  $L^0(\mathcal{F},B)$).

	\begin{cor} \label{cor.2.9}
		Let $( B , \| \cdot \|)$ be a Banach space over $\mathbb{K}$ and $V$ a
		nonempty weakly compact convex subset of $B$ such that $V$ has complete normal
		structure. Then, $L^0(\mathcal{F},V)$ is an $L^0$-convexly compact subset of
		$L^0(\mathcal{F},B)$ with random complete normal structure.
	\end{cor}
	
	\begin{proof}
		\cite[Corollary 1]{lim1974} shows that $V$ has normal structure, which,
		combined with \cite[Theorems 2.11 and 2.16]{GZWY20}, implies that
		$L^0(\mathcal{F},V)$ is an $L^0$-convexly compact subset with random normal
		structure. By Theorem \ref{thm.1.14}, $L^0(\mathcal{F},V)$ has random complete normal structure.
	\end{proof}

	\section{Proof of Theorem \ref{thm.1.15} }\label{sec.3}

	The aim of  this section is to prove Theorem \ref{thm.1.15}. To this end, we first present Lemma \ref{lem.3.2} below, 
	which is the special case of Theorem \ref{thm.1.15} when the family is finite.

	\par Lemmas \ref{lem.3.1} and \ref{lem.3.2}  below are respectively a random generalization of \cite[Theorem 4]{BK1966} and \cite[Theorem 3]{BK1966}, their proofs are omitted since the ideas of the proofs are respectively the same as \cite[Theorem 4]{BK1966} and \cite[Theorem 3]{BK1966} except for some changes in the random  setting as made in the proof of \cite[Theorem 3.8]{GZWG20}.	
	\begin{lem}\label{lem.3.1}
		Let $( E , \| \cdot \|)$ be a $\mathcal{T}_{\varepsilon,\lambda}$-complete
		$RN$ module over $\mathbb{K}$ with base $(\Omega,\mathcal{F},P)$, $G$
		an a.s. bounded $\mathcal{T}_{\varepsilon,\lambda}$-closed $L^0$-convex subset
		of $E$ such that $G$ has random normal structure, $M$ an $L^0$-convexly compact
		subset of $G$ and $T: G\rightarrow G$ a nonexpansive mapping with the property
		that  $[\{T^n(x): n=1,2,\cdots\} ]^{-}_{\varepsilon,\lambda}\cap M\neq
		\emptyset$ for any $x\in G$. Then $T$ has a fixed point.
	\end{lem}
	
	\begin{lem} \label{lem.3.2}
		Let $( E , \| \cdot \|)$ be a $\mathcal{T}_{\varepsilon,\lambda}$-complete
		$RN$ module over $\mathbb{K}$ with base $(\Omega,\mathcal{F},P)$, $G$
		an $L^0$-convexly compact subset with random normal structure of $E$ and
		$\mathcal{T}$ a finite commutative family of nonexpansive mappings from $G$ to
		$G$. Then $\mathcal{T}$ has a common fixed point.
	\end{lem}

	\par  It is obvious that a singleton is a $B_{\mathcal{F}}$-stable set in a natural way, whereas a nonempty set containing more than one point may be not necessarily a $B_{\mathcal{F}}$-stable set, Example \ref{ex.3.3} 
	below shows that we can always stabilize such a nonempty set by generating its $B_{\mathcal{F}}$-stable hull. 
	
	\begin{ex}\label{ex.3.3}
		Let $E$ be a nonempty set containing more than one point, denoting by $(e_{n},a_{n})_{n}$ a sequence
		$\{(e_{n},a_{n}), n\in \mathbb{N}\}$ in $E\times B_{\mathcal{F}}$, 
		and $S=\{(e_{n},a_{n})_{n}:\{a_{n},n\in\mathbb{N}\}$ is a partition of unity in $B_{\mathcal{F}}$, $\{e_{n},n\in
		\mathbb{N}\}$ is a sequence in $E\}$. Let us first define an equivalence relation
		$\sim_{1}$ on $S$ by $(e_{n},a_{n})_{n}\sim_{1}(f_{m},b_{m})_{m}$ iff
		$\bigvee\{a_{n}:e_{n}=z\}=\bigvee\{b_{m}:f_{m}=z\}$ for any $z\in E$. Denote by $[(e_{n},a_{n})_{n}]$ the
		equivalence class of $(e_{n},a_{n})_{n}$ under $\sim_{1}$ and assume $X=\{[(e_{n},a_{n})_{n}]:(e_{n},a_{n})_{n}\in
		S\}$, now let us define another equivalence relation $\sim_{2}$ on $X\times B_{\mathcal{F}}$
		by $([(e_{n},a_{n})_{n}],a)\sim_{2}([(f_{m},b_{m})_{m}],b)$ iff $a=b$
		and $(\bigvee\{a_{n}:e_{n}=z\})\wedge a=(\bigvee\{b_{m}:f_{m}=z\})\wedge a$ for
		any $z\in E$, then it is easy to check that $\sim_{2}$ is regular and $X$ is $B_{\mathcal{F}}$-stable
		with respect to $\sim_{2}$.  For any
		$x,y\in E$,  it is clear that $[x,1]=[y,1]$ iff $ x=y $, and hence $ [x,1] $ can be identified with $ x $ for any $x\in E$.
		 If each $(e_{n},a_{n})_{n}$ in $S$ can be interpreted as the step function taking the
		value $e_{n}\in E$ on $a_{n}$, then $ X $ can be interpreted as the set of equivalence classes of step functions under $\sim_{1}$.
		Besides, it is obvious that $[(e_{n},a_{n})_{n}]=\sum_{n=1}^{\infty}[e_{n},1]|a_{n}$, so we can use $
		\sum_{n=1}^{\infty}e_{n}|a_{n}$ for $ [(e_{n},a_{n})_{n}] $. We call $X$ the $B_{\mathcal{F}}$-stable hull
		of $E$,  denoted by $B_{\sigma}(E)$ in accordance with Definition \ref{defn.1.2}.	
	\end{ex}

	\par Although  a net $\{G_{\alpha},\alpha\in D\}$ of $\sigma$-stable subsets of an $L^{0}(\mathcal{F},\mathbb{K})$-module 
	is not necessarily a consistent net,  by the method provided by \cite{{GWT23}}, we can generate a consistent net from $\{G_{\alpha},\alpha\in D\}$. First, let $B_{\sigma}(D)$ be the
	$B_{\mathcal{F}}$-stable hull of a directed set $D$, then the directed relation $\leq$ on $D$ 
	naturally induces a directed relation on $B_{\sigma}(D)$ (also denoted by $\leq$) by
	$\sum_{n=1}^{\infty}\alpha_n|a_{n} \leq \sum_{n=1}^{\infty}\beta_m|b_{m}$ iff
	$\alpha_n \leq \beta_m$ whenever $a_n \wedge b_m>0$, it is known from \cite{{GWT23}} that $(B_{\sigma}(D),\leq)$
	is a $B_{\mathcal{F}}$-stable directed set. Now, for any  $\hat{\alpha}=\sum_{n=1}^{\infty}\alpha_{n}|a_{n}\in
	B_{\sigma}(D)$, where $\{\alpha_{n},n\in \mathbb{N}\}$  is a sequence in $D$ and $\{a_n,n\in \mathbb{N}\}$ a partition  
	of unity in $B_{\mathcal{F}}$, define $G_{\hat{\alpha}}=\sum_{n=1}^{\infty}\tilde{I}_{A_{n}}G_{\alpha_{n}}$, then 
	it is easy to check that $\{G_{\hat{\alpha}}, \hat{\alpha}\in
	B_{\sigma}(D)\}$ is a consistent net, which is called the consistent net generated by
	$\{G_{\alpha},\alpha\in D\}$.

	\par Lemmas \ref{lem.3.4} and \ref{lem.3.5} below are given to
	simplify the proof of Theorem \ref{thm.1.15}.

	\begin{lem}\label{lem.3.4}
		Let $( E, \| \cdot \|)$ be a $\mathcal{T}_{\varepsilon,\lambda}$-complete $RN$
		module over $\mathbb{K}$ with base $(\Omega,\mathcal{F},P)$, $K$ an
		$L^0$-convexly compact subset of $E$ and $\{M_{\alpha}, \alpha \in D\}$ a
		decreasing net of $\sigma$-stable subsets of $K$ with $R(M_{\alpha},K)=R(K,K)$ for any $
		\alpha \in D$. Then the consistent net $\{ M_{\hat{\alpha}}, \hat{\alpha} \in
		B_{\sigma}(D) \}$ generated by $\{M_{\alpha}, \alpha \in D\}$ is still  decreasing  and  satisfies the following conditions:
		\begin{enumerate}[(1)]
			\item $R(M_{\hat{\alpha}},K)=R(K,K)~\forall~ \hat{\alpha} \in
			B_{\sigma}(D)$.
			\item $\bigcup_{\hat{\alpha} \in
				B_{\sigma}(D)}C(M_{\hat{\alpha}},K) $ is $L^0$-convex.
		\end{enumerate}
	\end{lem}
	
	\begin{proof}
		For any $\hat{\alpha}=\sum_{n=1}^{\infty}\alpha_n|a_{n} $ and
		$\hat{\beta}=\sum_{m=1}^{\infty}\beta_m|b_{m} $ in $B_{\sigma}(D)$ with
		$\hat{\alpha} \leq \hat{\beta}$, where $\{a_{n},n\in \mathbb{N}\}$ and $\{b_{n},n\in \mathbb{N}\}$ are two partitions  of unity in
		$B_{\mathcal{F}}$ and $\{\alpha_n, n\in \mathbb{N}\}$  and  $\{\beta_n, n\in \mathbb{N}\}$ are 
		two sequences in $D$, we have 
		\begin{align}
			M_{\hat{\beta}}&=\sum_{m=1}^{\infty}\tilde{I}_{B_{m}}M_{\beta_{m}}\nonumber\\
			&=\sum_{m,n=1}^{\infty}\tilde{I}_{A_{n}\cap B_{m}}M_{\beta_{m}}\nonumber\\
			&\subset \sum_{m,n=1}^{\infty}\tilde{I}_{A_{n}\cap B_{m}}M_{\alpha_{n}}\nonumber\\
			&=M_{\hat{\alpha}},	\nonumber
		\end{align}
        which implies that $\{ M_{\hat{\alpha}}, \hat{\alpha} \in B_{\sigma}(D) \}$ is
		decreasing.
		
		(1). For any $\hat{\alpha} \in  B_{\sigma}(D)$, there exist a partition $\{a_{n},n\in \mathbb{N}\}$ of unity in
		$B_{\mathcal{F}}$ and a sequence $\{\alpha_n, n\in \mathbb{N}\}$ in $D$
		such that $\hat{\alpha}=\sum_{n=1}^{\infty}\alpha_n|a_{n}$, by (6) of
		Proposition \ref{prop.2.1}, we have
		\begin{align}
			R(M_{\hat{\alpha}},K)&=R(\sum_{n=1}^{\infty}\tilde{I}_{A_n}
			M_{\alpha_{n}},K)\nonumber\\
			&=\sum_{n=1}^{\infty}\tilde{I}_{A_n}R(M_{\alpha_{n}},K)\nonumber\\
			&=R(K,K).\nonumber
		\end{align}
		
		(2). For any $x,y \in \bigcup_{\hat{\alpha} \in
			B_{\sigma}(D)}C(M_{\hat{\alpha}},K)$, there exist $\hat{\alpha}_x,
		\hat{\alpha}_y \in B_{\sigma}(D)$ such that
		\begin{equation}\nonumber
			x\in C( M_{\hat{\alpha}_x}, K) ~\text{and}~ y\in C(
			M_{\hat{\alpha}_y},K).
		\end{equation}
		Let $\hat{\alpha}_x=\sum_{n=1}^{\infty}\alpha_n|a_{n}$ and $\hat{\alpha}_y=\sum_{m=1}^{\infty}\beta_m|b_{m}$, 
		where $\{a_{n},n\in \mathbb{N}\}$ and $\{b_{n},n\in \mathbb{N}\}$ are two partitions  of unity in
		$B_{\mathcal{F}}$ and $\{\alpha_n, n\in \mathbb{N}\}$  and  $\{\beta_n, n\in \mathbb{N}\}$ are 
		two sequences in $D$. By (7) of Proposition \ref{prop.2.1}, we have
		\begin{equation}\nonumber
			C( M_{\hat{\alpha}_x}, K)=C(
			\sum_{n=1}^{\infty}\tilde{I}_{A_n} M_{\alpha_n},
			K)=\sum_{n=1}^{\infty}\tilde{I}_{A_n} C( M_{\alpha_n}, K).
		\end{equation}
		and
		\begin{equation}\nonumber
			C(
			M_{\hat{\alpha}_y},K)=C(\sum_{m=1}^{\infty}\tilde{I}_{B_m}
			M_{\beta_m}, K)=\sum_{m=1}^{\infty}\tilde{I}_{B_m} C( M_{\beta_m},K).
		\end{equation}
		Since $D$ is a directed set, for each $n$ and $m$ such that $a_n\wedge b_m>0$,
		there exists $\gamma_{n,m} \in D$ such that $\gamma_{n,m}\geq \alpha_n$ and
		$\gamma_{n,m}\geq \beta_m$.
		Further, for any $u\in C( M_{\alpha_n},K)$  and any $v\in C( M_{\beta_m}, K)$, since $R(M_{\alpha},K)=R(K,K)$ for any $\alpha \in D$,  we have
		\begin{equation}\nonumber
			R(M_{\gamma_{n,m}},u)\leq
			R(M_{\alpha_n},u)=R(M_{\alpha_n},K)=R(K,K)=R(M_{\gamma_{n,m}},K)
		\end{equation}
		 and
		\begin{equation}\nonumber
			R(M_{\gamma_{n,m}},v)\leq
			R(M_{\beta_n},v)=R(M_{\beta_n},K)=R(K,K)=R(M_{\gamma_{n,m}},K),
		\end{equation}
		 which implies that
		$u\in C( M_{\gamma_{n,m}},K)~\text{and}~v\in C(
		M_{\gamma_{n,m}},K),$
		thus,
		\begin{equation}\nonumber
			C(  M_{\alpha_n},K) \subset C( M_{\gamma_{n,m}}, K)
			~\text{and}~ C( M_{\beta_m},K) \subset C( M_{\gamma_{n,m}},
			K).
		\end{equation}
		Let $\hat{\alpha}_{xy}=\sum_{n,m=1}^{\infty}\gamma_{n,m}|a_n \wedge b_m$, then $ \hat{\alpha}_{xy}\in
		B_{\sigma}(D)$ and we have
		\begin{align}\nonumber
			C( M_{\hat{\alpha}_x}, K)&=\sum_{n,m=1}^{\infty}\tilde{I}_{A_n \cap
				B_m} C(  M_{\alpha_n},K)\nonumber\\
			&\subset \sum_{n,m=1}^{\infty}\tilde{I}_{A_n \cap B_m} C(
			M_{\gamma_{n,m}}, K)\nonumber\\
			&=C( \sum_{n,m=1}^{\infty}\tilde{I}_{A_n \cap B_m}
			M_{\gamma_{n,m}},K)\nonumber \\
			&=C(  M_{\hat{\alpha}_{xy} },K) \nonumber
		\end{align}
		and
		\begin{align}\nonumber
			C( M_{\hat{\alpha}_y},K)&=\sum_{n,m=1}^{\infty}\tilde{I}_{A_n \cap
				B_m} C( M_{\beta_n}, K)\nonumber\\
			&\subset \sum_{n,m=1}^{\infty}\tilde{I}_{A_n \cap B_m} C(
			M_{\gamma_{n,m}}, K)\nonumber\\
			&=C( \sum_{n,m=1}^{\infty}\tilde{I}_{A_n \cap B_m}
			M_{\gamma_{n,m}},K)\nonumber \\
			&=C(  M_{\hat{\alpha}_{xy} },K). \nonumber
		\end{align}
		Since
		$C( M_{\hat{\alpha}_{xy}},K)$ is $L^0$-convex by Proposition
		\ref{prop.2.3}, then for any $\lambda \in
		L^0_{+}(\mathcal{F})$ with $0\leq \lambda \leq 1$, we have
		\begin{align}
			\lambda x+(1-\lambda)y &\in \lambda C(M_{\hat{\alpha}_x}, K)+(1-\lambda) C( M_{\hat{\alpha}_y},K)\nonumber\\
			&\subset  \lambda C( M_{\hat{\alpha}_{xy} },K)+(1-\lambda) C(M_{\hat{\alpha}_{xy} },K)\nonumber\\
			&= C( M_{\hat{\alpha}_{xy} },K)\nonumber\\
			&\subset \bigcup_{\hat{\alpha} \in B_{\sigma}(D)} C(M_{\hat{\alpha}},K),\nonumber
		\end{align}
		namely, $\bigcup_{\hat{\alpha} \in B_{\sigma}(D)} C(M_{\hat{\alpha}},K) $ is
		$L^0$-convex.
	\end{proof}

	\begin{lem}\label{lem.3.5}
		Let $( E , \| \cdot \|)$ be a $\mathcal{T}_{\varepsilon,\lambda}$-complete
		$RN$ module over $\mathbb{K}$ with base $(\Omega,\mathcal{F},P)$, $G$
		an $L^0$-convexly compact subset of $E$ and $\{M_{\alpha},\alpha \in D\}$ a net
		of nonempty $ \sigma $-stable subsets of $G$ with finite intersection property. Then
		\begin{equation}\nonumber
			\bigcap_{\hat{\alpha} \in B_{\sigma}(D)}  [Conv_{L^0} (M_{\hat{\alpha}})]^{-}_{\varepsilon,\lambda}\neq \emptyset,
		\end{equation}
		where $\{M_{\hat{\alpha}}, \hat{\alpha} \in B_{\sigma}(D)\}$ is the consistent net
		generated by $\{M_{\alpha},\alpha \in D\}$.
	\end{lem}
	\begin{proof}
		Let $B_\alpha=[Conv_{L^0} (M_{\alpha})]^{-}_{\varepsilon,\lambda}$ for each
		$\alpha \in D$ and $\{B_{\hat{\alpha}},\hat{\alpha} \in B_{\sigma}(D) \}$  the consistent net generated 
		by $\{ B_\alpha, \alpha \in D\}$.  It is easy to check that  $\tilde{I}_{A} H^{-}_{\varepsilon,\lambda}= (\tilde{I}_{A} H)
		^{-}_{\varepsilon,\lambda}$ and $\tilde{I}_{A} Conv_{L^{0}}(H)= Conv_{L^{0}}(\tilde{I}_{A} H)$ for any $A\in \mathcal{F}$ and any $H\subset E$, then for any $\hat{\alpha}=\sum_{n=1}^{\infty}\alpha_{n}|a_{n}$ in 
		$B_{\sigma}(D)$ with $\{a_{n},n\in \mathbb{N}\}$ a partition  of unity in $B_{\mathcal{F}}$ and 
		$\{\alpha_{n},n\in \mathbb{N}\}$ a sequence  in $D$, we have
		\begin{align}
			B_{\hat{\alpha}}&=\sum_{n=1}^{\infty}\tilde{I}_{A_{n}}B_{\alpha_{n}}\nonumber\\
			&=\sum_{n=1}^{\infty}\tilde{I}_{A_{n}}[Conv_{L^0} (M_{\alpha_{n}})]^{-}_{\varepsilon,\lambda}\nonumber\\
			&=\sum_{n=1}^{\infty}\tilde{I}_{A_{n}}[\tilde{I}_{A_{n}}Conv_{L^0} (\sum_{n=1}^{\infty}\tilde{I}_{A_{n}}M_{\alpha_{n}})]^{-}_{\varepsilon,\lambda}\nonumber\\
			&=\sum_{n=1}^{\infty}\tilde{I}_{A_{n}}[Conv_{L^0} (\sum_{n=1}^{\infty}\tilde{I}_{A_{n}}M_{\alpha_{n}})]^{-}_{\varepsilon,\lambda}\nonumber\\
			&=\sum_{n=1}^{\infty}\tilde{I}_{A_{n}}[Conv_{L^0} (M_{\hat{\alpha}})]^{-}_{\varepsilon,\lambda}\nonumber\\
			&=[Conv_{L^0} (M_{\hat{\alpha}})]^{-}_{\varepsilon,\lambda}.\nonumber
		\end{align}
		Then, we only need to prove $\bigcap_{\hat{\alpha} \in B_{\sigma}(D)}B_{\hat{\alpha}}\neq \emptyset$.
		
		Since $\{ B_\alpha, \alpha \in D\}$ is a net of nonempty
		$\mathcal{T}_{\varepsilon,\lambda}$-closed $L^0$-convex subsets of $G$ with
		finite intersection property,  we have 
		$\bigcap_{\alpha \in D}B_\alpha\neq \emptyset$ by the $L^0$-convex compactness of $G$. It is easy to check that $ \bigcap_{\alpha \in D} B_\alpha=\bigcap_{\hat{\alpha} \in B_{\sigma}(D)}B_{\hat{\alpha}}$, so the latter is nonempty.
	\end{proof}

	\par For an $RN$ module  $(E,\|\cdot\|)$ over $\mathbb{K}$ with base
	$(\Omega,\mathcal{F},P)$, $(E,\|\cdot\|)$ can also be endowed with the following locally $L^0$-convex topology, which will be used in the proof of Theorem \ref{thm.1.15}. For any $\varepsilon\in L^0_{++}(\mathcal{F})$, let
	$V(\theta, \varepsilon)=\{x\in E: \|x\|< \varepsilon ~\text{on}~ \Omega \}$,
	then $\mathcal{V}:=\{x+V(\theta, \varepsilon): x\in E \text{~and~}\varepsilon\in
	L^0_{++}(\mathcal{F}) \}$  forms a topological base of some Hausdorff topology for $ E $, called the locally $L^0$-convex topology \cite{FKV2009},
	denoted by $\mathcal{T}_{c}$. It is known from \cite{FKV2009} that $ (E,\mathcal{T}_{c}) $ is a topological module over the topological ring $(L^{0}(\mathcal{F},\mathbb{K}),\mathcal{T}_{c})$. 
	As shown by \cite[Theorem 3.12]{Guo10}, if $G$ is a $\sigma$-stable subset of an $RN$ module, then
	$G^{-}_{\varepsilon,\lambda}=G^{-}_{c}$, where
	$G^{-}_{\varepsilon,\lambda}$ and $G^{-}_{c}$ denote the closure of $G$
	under $\mathcal{T}_{\varepsilon, \lambda}$ and $\mathcal{T}_{c}$, respectively.
	
	\par Let $(E,\|\cdot\|)$ be an $RN$ module 
	with base $(\Omega,\mathcal{F},P)$, $G$ and $F$ be two nonempty $\sigma$-stable subsets of
	$E$. The mapping $T: G\rightarrow F$ is said to be $\sigma$-stable if 
	$T(\sum_{n=1}^{\infty}\tilde{I}_{A_n}x_n)=\sum_{n=1}^{\infty}\tilde{I}_{A_n}T(x_n)$
	for each sequence $\{x_n, n\in \mathbb{N}\}$ in $G$ and each $\{A_n, n\in
	\mathbb{N}\} \in \Pi_{\mathcal{F}}$. It is known from \cite[Lemma 2.11]{GZWG20} that 
	if $E$ is $\mathcal{T}_{\varepsilon, \lambda}$-complete and $T$ is nonexpansive, then $T$ must be $\sigma$-stable.

	\par Now we are ready to prove Theorem \ref{thm.1.15}.

	\begin{proof}[\textbf{Proof of Theorem \ref{thm.1.15}}]
		Let $\mathcal{A}=\{H: H ~$is a nonempty$~
		\mathcal{T}_{\varepsilon,\lambda}$-closed$~L^0$-convex$ ~$subset of G such
		that$~ T(H)\subset H ~$for any$~ T\in  \mathcal{T}\}$, then $\mathcal{A}$ is a
		partially ordered set under the usual inclusion relation. By the $L^0$-convex
		compactness of $G$ and the Zorn's Lemma, there exists a minimal element $K$ in
		$\mathcal{A}$.  We assert that $K$ consists of a single point, which is the
		desired common fixed point of $\mathcal{T}$.
		
		\par
		Otherwise, $K$ contains more than one point. Let $\mathcal{Q}$ be the family of nonempty finite subsets of $\mathcal{T}$.
		For any $Q\in \mathcal{Q}$, let
		$$M_{Q}=\{x\in K: T(x)=x ~\forall~ T\in Q\},$$
		then  $M_{Q}$ is nonempty by Lemma \ref{lem.3.2} and is
		$\sigma$-stable since each $T$ is $\sigma$-stable. Obviously, $M_{Q_{2}}\subset M_{Q_{1}}$ if
		$Q_{1}\subset Q_{2}$.
		Therefore, $\{M_{Q},Q\in \mathcal{Q}\}$ is a decreasing net of nonempty
		$\sigma$-stable subsets of $K$, and then the consistent net $\{M_{\hat{Q}}, \hat{Q}\in
		\sigma(\mathcal{Q} )\}$ generated by $\{M_{Q},Q\in \mathcal{Q}\}$, is a
		decreasing consistent net of nonempty $\sigma$-stable subsets of $K$.  Since $G$ has random normal
		structure, by Theorem \ref{thm.1.14}, $G$ has random complete normal structure.
		The remainder of the proof is divided into two steps.
		
		\textbf{Step 1.} We  prove that $R(M_{\hat{Q}},K)=R(K,K)$ for any
		$\hat{Q}\in B_{\sigma}(\mathcal{Q})$.
		
		By Lemma \ref{lem.3.4}, it suffices to prove $R( M_{Q},K)=R( K,K)$ for any
		$Q\in\mathcal{Q}$.
		For any given $Q_0 \in \mathcal{Q}$, let $r_0:=R( M_{Q_0},K)$ and
		\begin{equation}\nonumber
			H_{Q}=\{x\in K: R( M_{Q},x)\leq r_0\} ~\text{for any}~ Q\in \mathcal{Q}
			~\text{with}~ Q_0\subset Q.
		\end{equation}
		Then $H_{Q_0}=C(M_{Q_0},K)\neq \emptyset$ since $K$ is $L^{0}$-convexly compact, and the collection
		$\{H_{Q}:Q\in \mathcal{Q} ~\text{with}~ Q_0\subset Q\}$ forms an increasing net
		of nonempty $L^0$-convex sets. Let
		$S_{Q_{0}}=\bigcup\{H_{Q}:Q\in \mathcal{Q} ~\text{and}~ Q_0\subset Q\},$ we
		claim that $[\sigma(S_{Q_{0}})]^{-}_{\varepsilon,\lambda}\in \mathcal{A}$.
		
		First, since each $H_{Q}$ is $L^{0}$-convex and $ H_{Q} $ increases with $Q$, it is clear that $S_{Q_{0}}$ is $L^{0}$-convex, and so is
		$\sigma(S_{Q_{0}})$, which implies that
		$[\sigma(S_{Q_{0}})]^{-}_{\varepsilon,\lambda}$ is a nonempty
		$\mathcal{T}_{\varepsilon,\lambda}$-closed $L^{0}$-convex subset of $G$.
		
		Second, for any $T\in \mathcal{T}$ and any $x\in S_{Q_{0}}$, there exists
		$Q\in \mathcal{Q}$ with $Q_0\subset Q$ such that $x\in H_{Q}$ and $T\in Q$,
		therefore,
		\begin{align}\nonumber
			R(M_{Q},T(x))&=\bigvee \{\|T(x)-z\|:z\in M_{Q}\} \nonumber \\
			&=\bigvee \{\|T(x)-T(z)\|:z\in M_{Q}\} \nonumber \\
			&\leq\bigvee \{\|x-z\|:z\in M_{Q}\} \nonumber \\
			&=R(M_{Q},x)\nonumber \\
			&\leq r_0, \nonumber
		\end{align}
		which implies that $T(x)\in H_{Q}\subset S_{Q_{0}}$. Thus $T(
		S_{Q_{0}})\subset S_{Q_{0}} $ for each $ T\in \mathcal{T}$. 
		Further, for each $T\in \mathcal{T}$, since $T$ is $\sigma$-stable, then $T(\sigma(
		S_{Q_{0}}))\subset \sigma( S_{Q_{0}})$, and hence $T([\sigma(
		S_{Q_{0}})]^{-}_{\varepsilon,\lambda})\subset [\sigma(
		S_{Q_{0}})]^{-}_{\varepsilon,\lambda}$ by the $ \mathcal{T}_{\varepsilon,\lambda} $-continuity of $T$, which means 
		$[\sigma(S_{Q_{0}})]^{-}_{\varepsilon,\lambda}=K$ by the minimality of $K$.

		Since $[\sigma( S_{Q_{0}})]^{-}_{c}=[\sigma(
		S_{Q_{0}})]^{-}_{\varepsilon,\lambda}=K$, for any $\varepsilon \in
		L^{0}_{++}(\mathcal{F})$ and any $x\in K$, there exists $h\in \sigma(
		S_{Q_{0}})$ such that $\|x-h\|\leq \varepsilon$. For such a point $h$, there exist
		$\{A_n,n\in \mathbb{N}\} \in \Pi_{\mathcal{F}}$ and $\{h_n,n\in N\}$ in
		$S_{Q_{0}}$ such that $h=\sum_{n=1}^{\infty}\tilde{I}_{A_n}h_n$.  For each
		$h_n$, there exists $Q_n \in \mathcal{Q}$ with $Q_0\subset Q_n$ such that
		$h_n\in H_{Q_n}$, namely, $R(M_{Q_n},h_n)\leq r_0$. Let
		$\hat{Q}_{x}=\sum_{n=1}^{\infty}Q_n|a_n$, then $\hat{Q}_{x}\in B_{\sigma}(\mathcal{Q})$ and  $Q_{0}\subset \hat{Q}_{x}$, we have
		\begin{align}\nonumber
			R(M_{\hat{Q}_x},x)&=\bigvee \{\|x-y\|:y\in M_{\hat{Q}_{x}}\} \nonumber \\
			&\leq \bigvee \{\|x-h\|+\|h-y\|: y\in M_{\hat{Q}_{x}}\}\nonumber \\
			&\leq \bigvee \{\|h-y\|:y\in M_{\hat{Q}_{x}}\} +\varepsilon \nonumber \\
			&=R(M_{\hat{Q}_{x}},h)+\varepsilon \nonumber \\
			&=R(
			\sum_{n=1}^{\infty}\tilde{I}_{A_n}M_{Q_n},\sum_{n=1}^{\infty}\tilde{I}_{A_n}
			h_n)+\varepsilon \nonumber\\
			&=\sum_{n=1}^{\infty}\tilde{I}_{A_n} R( M_{Q_n},h_n )+\varepsilon \nonumber
			\\
			&\leq r_0 +\varepsilon, \nonumber
		\end{align}
		which implies that $M_{\hat{Q}_x}\subset B(x,r_0+\varepsilon)$, where
		$B(x,r_0+\varepsilon):=\{y\in K: \|x-y\|\leq r_0+\varepsilon\}$ is
		$\mathcal{T}_{\varepsilon,\lambda}$-closed and $L^0$-convex. Thus, we have
		\begin{align}\nonumber
			&\bigcap_{\hat{Q}\in B_{\sigma}(\mathcal{Q}), Q_0
			\subset \hat{Q}}
			[Conv_{L^0}(M_{\hat{Q}})]^{-}_{\varepsilon,\lambda}\nonumber\\
		\subset	& \bigcap_{\hat{Q}_x\in B_{\sigma}(\mathcal{Q}), Q_0\subset \hat{Q}_x}
			[Conv_{L^0}(M_{\hat{Q}_x})]^{-}_{\varepsilon,\lambda} \nonumber \\
		\subset &\bigcap_{x\in K} B(x,r_0+\varepsilon).\nonumber
		\end{align}
		Further, Lemma \ref{lem.3.5} shows that
		$$ \bigcap_{\hat{Q}\subset B_{\sigma}(Q), Q_0\subset \hat{Q}}
		[Conv_{L^0}(M_{\hat{Q}})]^{-}_{\varepsilon,\lambda} \neq \emptyset.$$
		Then,  $\bigcap_{x\in K} B(x,r_0+\varepsilon)\neq \emptyset $. For any $y\in
		\bigcap_{x\in K} B(x,r_0+\varepsilon)$, we have $R(K,y)\leq r_0+\varepsilon$,
		which implies that $R(K,K)\leq R(K,y)\leq r_0+\varepsilon$. By the arbitrariness
		of $\varepsilon$, we have $R(K,K) \leq r_0=R( M_{Q_0},K)$. On the other hand,
		since $M_{Q_0}\subset K$,  we have $R( M_{Q_0},x)\leq R(K,x)$ for any $x\in K$,
		which implies that $R( M_{Q_0},K)\leq R(K,K)$.
		
		To sum up, $R( M_{Q},K)= R(K,K)$ for each $Q\in \mathcal{Q}$.

		\textbf{Step 2.} We prove that $[\bigcup_{\hat{Q}\in B_{\sigma}(\mathcal{Q})}
		C(  M_{\hat{Q}},K)]^{-}_{\varepsilon,\lambda}=K$, which leads to a contradiction 
		with the fact that $G$ has random complete normal structure.
		
		By the minimality of $K$, it suffices to prove $[\bigcup_{\hat{Q}\in
			B_{\sigma}(\mathcal{Q})} C(  M_{\hat{Q}},K)]^{-}_{\varepsilon,\lambda}\in
		\mathcal{A} $.
		
		First, Lemma \ref{lem.3.4} implies that $\bigcup_{\hat{Q}\in
			B_{\sigma}(\mathcal{Q})} C( M_{\hat{Q}}, K)$ is nonempty and $L^0$-convex,
		so $[\bigcup_{\hat{Q}\in B_{\sigma}(\mathcal{Q})} C(
		M_{\hat{Q}},K)]^{-}_{\varepsilon,\lambda} $ is a nonempty
		$\mathcal{T}_{\varepsilon,\lambda}$-closed $L^0$-convex subset of $K$.
		
		Second, we prove that
		$$T([\bigcup_{\hat{Q}\in B_{\sigma}(\mathcal{Q})}
		C(  M_{\hat{Q}},K)]^{-}_{\varepsilon,\lambda})\subset
		[\bigcup_{\hat{Q}\in B_{\sigma}(\mathcal{Q})} C(
		M_{\hat{Q}},K)]^{-}_{\varepsilon,\lambda}$$
		for any $T\in \mathcal{T}$.
		Since each $T$ is continuous, we only need to show that $T(\bigcup_{\hat{Q}\in
			B_{\sigma}(\mathcal{Q})} C( M_{\hat{Q}}, K))\subset \bigcup_{\hat{Q}\in
			B_{\sigma}(\mathcal{Q})} C( M_{\hat{Q}}, K)$ for any $T\in \mathcal{T}$.
		
		For any $T\in \mathcal{T}$ and any $x\in \bigcup_{\hat{Q}\in B_{\sigma}(\mathcal{Q})} C(
		M_{\hat{Q}}, K)$, there exists $\hat{Q}_{x}\in
		B_{\sigma}(\mathcal{Q})$ such that $x \in C( M_{\hat{Q}_x},K)$. For
		$\hat{Q}_{x}$, there exist $\{A_n, n\in \mathbb{N}\} \in \Pi_{\mathcal{F}}$ and
		$\{Q_n,n\in \mathbb{N}\}$ in $\mathcal{Q}$ such that
		$\hat{Q}_x=\sum_{n=1}^{\infty}Q_n|a_{n}$. Let $Q'_n=Q_n \cup \{T\}$ for each
		$n\in \mathbb{N}$ and $\hat{Q}'_x=\sum_{n=1}^{\infty}Q'_n|a_{n}$, then
		$\hat{Q}'_{x}\in B_{\sigma}(\mathcal{Q})$ and we have
		\begin{align}\nonumber
			C(M_{\hat{Q}_x},K)&=\sum_{n=1}^{\infty}\tilde{I}_{A_n}
			C(M_{Q_n}, K)\nonumber\\
			&\subset \sum_{n=1}^{\infty}\tilde{I}_{A_n} C( M_{Q'_n}, K)\nonumber\\
			&=C(\sum_{n=1}^{\infty}\tilde{I}_{A_n} M_{Q'_n}, K)\nonumber\\
			&=C(M_{\hat{Q}'_{x}}, K),\nonumber
		\end{align}
		which implies that $x\in C( M_{\hat{Q}'_{x}}, K)$. Further, since
		\begin{align}
			R( M_{Q'_n},T(x))&=\bigvee \{\|T(x)-z\|:z\in M_{Q'_n}\}\nonumber\\
			&=\bigvee \{\|T(x)-T(z)\|:z\in M_{Q'_n}\}\nonumber\\
			&\leq \bigvee \{\|x-z\|:z\in M_{Q'_n}\}\nonumber\\
			&= R( M_{Q'_n},x) \nonumber
		\end{align}
		for each $n\in \mathbb{N}$, then we have
		\begin{align}\nonumber
			R( M_{\hat{Q}'_{x}},T(x))&= R( \sum_{n=1}^{\infty}\tilde{I}_{A_n} M_{
				Q'_n},T(x))\nonumber\\
			&=\sum_{n=1}^{\infty}\tilde{I}_{A_n} R( M_{ Q'_n},T(x)) \nonumber\\
			&\leq \sum_{n=1}^{\infty}\tilde{I}_{A_n} R( M_{ Q'_n},x) \nonumber\\
			&= R( \sum_{n=1}^{\infty}\tilde{I}_{A_n}M_{Q'_n}, x)\nonumber\\
			&=R(  M_{\hat{Q}'_{x}},K),\nonumber
		\end{align}
		that is, $T(x) \in C( M_{\hat{Q}'_{x}}, K) \subset
		\bigcup_{\hat{Q}\in B_{\sigma}(\mathcal{Q})} C( M_{\hat{Q}}, K)$. Thus,
		$$T(\bigcup_{\hat{Q}\in B_{\sigma}(\mathcal{Q})} C( M_{\hat{Q}}, K))
		\subset \bigcup_{\hat{Q}\in B_{\sigma}(\mathcal{Q})} C( M_{\hat{Q}}, K)$$
		for each $T\in \mathcal{T}$.
	\end{proof}

	\par Corollary \ref{cor.3.6} below is a random generalization of \cite[Theorem 2]{Browder1965}.

	\begin{cor}\label{cor.3.6}
		Let $( E , \| \cdot \|)$ be a $\mathcal{T}_{\varepsilon,\lambda}$-complete
		random uniformly convex $RN$ module over $\mathbb{K}$ with base
		$(\Omega,\mathcal{F},P)$, $G$ a nonempty a.s. bounded
		$\mathcal{T}_{\varepsilon,\lambda}$-closed $L^0$-convex subset of $E$ and
		$\mathcal{T}$ a commutative family of nonexpansive mappings from $G$ to
		$G$. Then $\mathcal{T}$ has a common fixed point.
	\end{cor}

	\par Corollary \ref{cor.3.7} below is a random generalization of \cite[Theorem 1.9]{Raj2015}.

	\begin{cor}\label{cor.3.7}
		Let $( E , \| \cdot \|)$ be a $\mathcal{T}_{\varepsilon,\lambda}$-complete
		$RN$ module over $\mathbb{K}$ with base $(\Omega,\mathcal{F},P)$, $G$ an
		$L^0$-convexly compact subset with random normal structure of $E$ and
		$\mathcal{T}$ a commutative family of nonexpansive mappings from $G$ to $G$
		with $[T(G)]^{-}_{\varepsilon,\lambda}=G$ for each $T\in
		\mathcal{T}$. Then $\mathcal{T}$ has a common fixed point in $C(G,G)$.
	\end{cor}
	\begin{proof}
		Obviously, $G$ is a.s. bounded. For any $x\in C(G,G)$ and any $T\in
		\mathcal{T}$, we have
		$$R(G,Tx)=R([T(G)]^{-}_{\varepsilon,\lambda},Tx)=R(T(G),Tx)\leq
		R(G,x)=R(G,G),$$
		which implies that $Tx\in C(G,G)$, and then
		$T(C(G,G))\subset C(G,G)$ for each $T\in \mathcal{T}$.
		Further, Proposition \ref{prop.2.3} shows that $C(G,G)$ is an
		$L^0$-convexly compact subset with random normal structure. Thus,
		by Theorem
		\ref{thm.1.15}, there exists $x\in C(G,G)$ such that $Tx=x$ for each $T\in
		\mathcal{T}$.
	\end{proof}

	\section{Proofs of Theorems \ref{thm.1.16} and \ref{thm.1.17}}\label{sec.4}

	This section is devoted to the proofs of Theorems \ref{thm.1.16} and
	\ref{thm.1.17}, which are the random generalizations of the
	classical results in  \cite{Lim2003} and \cite{BM1948}, respectively.

	\par An $RN$ module $( E , \| \cdot \|)$ is said to be
	random strictly convex \cite{GZ10} if for any $x$ and $y$ in $E\backslash
	\{\theta \}$ such that $\|x+y\|=\|x\|+\|y\|$, there exists $\xi \in
	L^0_{+}(\mathcal{F})$ such that $\xi >0$ on  $A_{xy}$ and
	$\tilde{I}_{A_{xy}}x=\xi( \tilde{I}_{A_{xy}} y )$.

	\par The study of affine geometry in regular $L^0$-module began with \cite{WGL22}.  We give Propositions \ref{prop.4.1} and \ref{prop.4.2} below 
	since they are of independent interest. Proposition \ref{prop.4.1} is a basic result concerning $L^0$-affine mappings, which is a random generalization of \cite[Lemma 2.1]{Raj2015}.

	\begin{prop} \label{prop.4.1}
		Let $( E , \| \cdot \|)$ be a $\mathcal{T}_{\varepsilon,\lambda}$-complete
		random strictly convex $RN$ module over $\mathbb{K}$ with base
		$(\Omega,\mathcal{F},P)$, $G$ an $L^0$-convex subset of $E$ and $T:
		G\rightarrow E$ an isometric mapping, then $T$ is $L^0$-affine (namely, 
		$T(\lambda x+(1-\lambda)y)=\lambda T(x)+(1-\lambda) T(y)$ for any $x,y\in G$ and any $\lambda \in
		L^0_{+}(\mathcal{F})$ with $0\leq \lambda \leq 1$).
	\end{prop}
	
	\begin{proof}
		For any $x,y\in G$ and any $\alpha \in L^0_{+}(\mathcal{F})$ with $0\leq
		\alpha \leq 1$, let $z=\alpha x+(1-\alpha)y$, it suffices to prove that
		$Tz=\alpha Tx+(1-\alpha)Ty$.
		
		Since	
		$\|x-z\|=(1-\alpha)\|x-y\| $ and $\|z-y\|=\alpha\|x-y\|$,  it follows that
		\begin{align}\nonumber
			\|Tx-Tz\|+\|Tz-Ty\|&=\|x-z\|+\|z-y\|\nonumber\\
			&=\|x-y\|.\nonumber
		\end{align}
		Then $\|x-y\|>0$ on $A_{(Tx-Tz)(Tz-Ty)}$. By the definition of random strict
		convexity, there exists $\xi \in L^0_{+}(\mathcal{F})$ such that $\xi >0$ on
		$A_{(Tx-Tz)(Tz-Ty)}$ and
		$\tilde{I}_{A_{(Tx-Tz)(Tz-Ty)}}(Tx-Tz)=\xi \tilde{I}_{A_{(Tx-Tz)(Tz-Ty)}}
		(Tz-Ty)$, namely,
		$$\tilde{I}_{A_{(Tx-Tz)(Tz-Ty)}}Tz \\
		=\tilde{I}_{A_{(Tx-Tz)(Tz-Ty)}}(\frac{1}{1+\xi}Tx+\frac{\xi}{1+\xi}Ty).$$
		Hence,
		we have
		\begin{align}\nonumber
			&\tilde{I}_{A_{(Tx-Tz)(Tz-Ty)}}\|Tx-Tz\|-\xi \tilde{I}_{A_{(Tx-Tz)(Tz-Ty)}}
			\|Tz-Ty\| \nonumber\\
			=& \tilde{I}_{A_{(Tx-Tz)(Tz-Ty)}}\|x-z\|- \tilde{I}_{A_{(Tx-Tz)(Tz-Ty)}}\xi
			\|z-y\| \nonumber\\
			=& \tilde{I}_{A_{(Tx-Tz)(Tz-Ty)}}(1-\alpha)\|x-y\|-
			\tilde{I}_{A_{(Tx-Tz)(Tz-Ty)}}\xi \alpha\|x-y\|\nonumber\\
			=& \tilde{I}_{A_{(Tx-Tz)(Tz-Ty)}} (1-(\xi+1) \alpha) \|x-y\|\nonumber\\
			=& 0,\nonumber
		\end{align}
		which implies that  $\tilde{I}_{A_{(Tx-Tz)(Tz-Ty)}}
		\frac{1}{1+\xi}=\tilde{I}_{A_{(Tx-Tz)(Tz-Ty)}}\alpha$. Therefore,
		\begin{equation}\nonumber
			\tilde{I}_{A_{(Tx-Tz)(Tz-Ty)}}Tz=\tilde{I}_{A_{(Tx-Tz)(Tz-Ty)}}(\alpha
			Tx+(1-\alpha)Ty).
		\end{equation}\nonumber
		
		On the other hand,
		\begin{align}\nonumber
			\tilde{I}_{A^c_{(Tx-Tz)}\cap
				A_{(Tz-Ty)}}\|Tx-Tz\|&=\tilde{I}_{A^c_{(Tx-Tz)}\cap
				A_{(Tz-Ty)}}\|x-z\|\nonumber\\
			&=\tilde{I}_{A^c_{(Tx-Tz)}\cap A_{(Tz-Ty)}} (1-\alpha) \|x-y\|\nonumber\\
			&=0\nonumber
		\end{align}
		implies that $\tilde{I}_{A^c_{(Tx-Tz)}\cap A_{(Tz-Ty)}} (1-\alpha)=0$. Then,
		we have
		\begin{align}\nonumber
			\tilde{I}_{A^c_{(Tx-Tz)}\cap A_{(Tz-Ty)}} Tz&=\tilde{I}_{A^c_{(Tx-Tz)}\cap
				A_{(Tz-Ty)}} Tx \nonumber\\
			&=\tilde{I}_{A^c_{(Tx-Tz)}\cap A_{(Tz-Ty)}} (\alpha
			Tx+(1-\alpha)Ty).\nonumber
		\end{align}
		
		Similarly, we can obtain $$\tilde{I}_{A_{(Tx-Tz)}\cap A^c_{(Tz-Ty)}}
		Tz=\tilde{I}_{A_{(Tx-Tz)}\cap A^c_{(Tz-Ty)}} (\alpha Tx+(1-\alpha)Ty)$$ and
		$$\tilde{I}_{A^c_{(Tx-Tz)}\cap A^c_{(Tz-Ty)}} Tz=\tilde{I}_{A^c_{(Tx-Tz)}\cap
			A^c_{(Tz-Ty)}} (\alpha Tx+(1-\alpha)Ty).$$
		
		Since $\{A_{(Tx-Tz)(Tz-Ty)},A^c_{(Tx-Tz)}\cap A_{(Tz-Ty)},A_{(Tx-Tz)}\cap
		A^c_{(Tz-Ty)},\\A^c_{(Tx-Tz)}\cap A^c_{(Tz-Ty)}\}  \in \Pi_{\mathcal{F}}$, we have
		$Tz=\alpha Tx+(1-\alpha)Ty$.
	\end{proof}

	\begin{prop}\label{prop.4.2}
		Let $( E , \| \cdot \|)$ be a $\mathcal{T}_{\varepsilon,\lambda}$-complete
		random uniformly convex $RN$ module over $\mathbb{K}$ with base
		$(\Omega,\mathcal{F},P)$, $G$ a nonempty a.s. bounded
		$\mathcal{T}_{\varepsilon,\lambda}$-closed $L^0$-convex subset of $E$ and $T:
		G\rightarrow G$ an isometric mapping, then
		$R(T(G),T(G))=R(G,G)$ and $C(T(G),T(G))=T( C(G,G) )$.
	\end{prop}
	\begin{proof}
		It is  known from \cite[Proposition 4.4]{GZ10} that a random uniformly convex $RN$ module 
		must be random strictly convex, then $T$ is $L^0$-affine by Proposition \ref{prop.4.1}, and hence  $T(G)$ is $L^0$-convex.
		It is clear that $T(G)$ is an a.s. bounded
		$\mathcal{T}_{\varepsilon,\lambda}$-closed $L^0$-convex subset of $G$, as in the proof of Corollary \ref{cor.2.8}, 
		$T(G)$ is $L^0$-convexly compact, and further $C(T(G),T(G))\neq
		\emptyset$ by Proposition \ref{prop.2.3}. Then, we have
		\begin{align}\nonumber
			R(T(G),T(G))&=\bigwedge \{R(T(G),Tx):x\in G\}\nonumber\\
			&=\bigwedge\{R(G,x):x\in G\} \nonumber\\
			&=R(G,G).\nonumber	
		\end{align}
		
		For any $x_{0}\in C(G,G)$, since
		$$R(T(G),Tx_{0})=R(G,x_0)=R(G,G)=R(T(G),T(G)),$$
		 then $Tx_{0}\in
		C(T(G),T(G))$, which implies that $T( C(G,G) )\subset
		C(T(G),T(G))$.
		Conversely, for any $y' \in C(T(G),T(G))$, there exists $x'\in G$
		such that $y'=Tx'$ and 
		$$R(G,x')=R(T(G),y')= R(T(G),T(G))=R(G,G),$$
		 then $x'\in
		C(G,G)$, which implies that $C(T(G),T(G)) \subset T(
		C(G,G) )$.
		Thus, $C(T(G),T(G))=T(C(G,G) )$.
	\end{proof}

	\par We first give Lemmas \ref{lem.4.3} and \ref{lem.4.4} below for the proofs of Lemma \ref{lem.4.5} and Theorem \ref{thm.1.16}.

	\begin{lem}\label{lem.4.3}
	Let $E$, $G$ and $\mathcal{E}$ be the same as in Definition \ref{defn.2.5},
	we have the following statements:
	\begin{enumerate}[(1)]
	\item $AR(\mathcal{E},\cdot)$ is $\mathcal{T}_{\varepsilon,\lambda}$-continuous and $L^0$-convex.
	\item  $\tilde{I}_{A} AR(\mathcal{E},x)=AR(\tilde{I}_{A}\mathcal{E},\tilde{I}_{A}x)$ for any $x\in
			E$ and any $A\in \mathcal{F}$. Further, if $E$ is $\sigma$-stable, then
			$$AR(\mathcal{E},\sum_{n=1}^{\infty}\tilde{I}_{A_n}
			x_n)=\sum_{n=1}^{\infty}\tilde{I}_{A_n} AR(\mathcal{E},x_n)$$
			 for any sequence $\{x_n,n\in
			\mathbb{N}\}$ in $E$ and any $\{A_{n},n\in \mathbb{N}\} \in
			\Pi_{\mathcal{F}}$.
	\item $\tilde{I}_{A} AR(\mathcal{E},G)=AR(\tilde{I}_{A}\mathcal{E},\tilde{I}_{A} G)$ for any $A\in
			\mathcal{F}$. Further, if $E$ is $\sigma$-stable, then
			$$AR(\mathcal{E},
			\sum_{n=1}^{\infty}\tilde{I}_{A_n}G_{n})=\sum_{n=1}^{\infty}\tilde{I}_{A_n} AR(\mathcal{E},G_{n})$$ 
			for any sequence $\{G_{n},n\in
			\mathbb{N}\}$ of  nonempty subsets of $E$  and any $\{A_{n},n\in
			\mathbb{N}\} \in \Pi_{\mathcal{F}}$.
	\item If $G$ is finitely stable and $AC(H,G)\neq \emptyset$, then
			$\tilde{I}_{A} AC(\mathcal{E},G)= AC(\tilde{I}_{A}\mathcal{E},\tilde{I}_{A}G)$ for
			any $A\in \mathcal{F}$. Further, if $E$ is $\sigma$-stable, then 
			$$AC(\mathcal{E}, \sum_{n=1}^{\infty}\tilde{I}_{A_n}G_{n})=\sum_{n=1}^{\infty}\tilde{I}_{A_n}
			AC(\mathcal{E},G_{n})$$
			 for any sequence $\{G_{n},n\in
			\mathbb{N}\}$ of  finitely stable subsets of $E$ and any $\{A_{n},n\in \mathbb{N}\} \in
			\Pi_{\mathcal{F}}$.
		\end{enumerate}
	\end{lem}
	\begin{proof}
		By the same methods used in the proof of Proposition \ref{prop.2.1}, it is easy to check that 
		$(2)$, $(3)$ and $(4)$ are true. We focus on proving (1).
		
		By simple calculations, we have
		$|AR(\mathcal{E},x)-AR(\mathcal{E},y)|\leq \|x-y\|$ for any $x,y\in W$,
		thus $AR(\mathcal{E},\cdot)$ is
		$\mathcal{T}_{\varepsilon,\lambda}$-continuous. For any $\alpha_1, \alpha_2 \in
		\Lambda$ and any $\lambda \in L^0_{+}(\mathcal{F})$ with $0\leq \lambda\leq1$,
		there exists $\alpha_3 \in \Lambda$ with $\alpha_3\geq \alpha_1$ and
		$\alpha_3\geq \alpha_2$, such that $W_{\alpha_3}\subset \lambda
		W_{\alpha_1}+(1-\lambda)W_{\alpha_2}$, we have
		\begin{align}\nonumber
			R(W_{\alpha_3},\lambda x+(1-\lambda)y)&\leq R(\lambda
			W_{\alpha_1}+(1-\lambda)W_{\alpha_2},\lambda x+(1-\lambda)y)\\ \nonumber
			&=\bigvee \{ \|\lambda x+(1-\lambda)y-\lambda a-(1-\lambda)b\|: a\in
			W_{\alpha_1}, b\in W_{\alpha_2}\} \\ \nonumber
			&\leq \lambda R(W_{\alpha_1},x)+(1-\lambda) R(W_{\alpha_2},y). \nonumber
		\end{align}
		From 
		\begin{align}\nonumber
			AR( \mathcal{E},\lambda x+(1-\lambda)y)&\leq
			R(W_{\alpha_3},\lambda x+(1-\lambda)y)\\ \nonumber
			&\leq  \lambda  R(W_{\alpha_1},x)+
			(1-\lambda) R(W_{\alpha_2},y),  \nonumber
		\end{align} 
		we have $AR( \mathcal{E},\lambda x+(1-\lambda)y)\leq \lambda AR(\mathcal{E},x) + (1-\lambda)AR(\mathcal{E},y) $ by the arbitrariness of $\alpha_1$ and $\alpha_2$.
	\end{proof}
	
	\par Similarly to Proposition \ref{prop.2.3}, one can have the following. 
	
	\begin{lem}\label{lem.4.4}
		Let $( E , \| \cdot \|)$ be a $\mathcal{T}_{\varepsilon,\lambda}$-complete
		$RN$ module over $\mathbb{K}$ with base $(\Omega,\mathcal{F},P)$, $G$
		an $L^{0}$-convexly compact subset of $E$ and  $\mathcal{E}=\{B_\alpha, \alpha\in
		\Lambda\}$ a decreasing net of a.s. nonempty bounded subsets of $E$. Then $AC(\mathcal{E},G)$ is a  nonempty
		$\mathcal{T}_{\varepsilon,\lambda}$-closed $L^0$-convex subset of $G$.
	\end{lem}
	
	\par Lemma \ref{lem.4.5} below is crucial for the proof of Theorem \ref{thm.1.16}.
	
	\begin{lem} \label{lem.4.5}
		Let $( E , \| \cdot \|)$ be a $\mathcal{T}_{\varepsilon,\lambda}$-complete
		$RN$ module over $\mathbb{K}$ with base $(\Omega,\mathcal{F},P)$, $G$
		an $L^0$-convexly compact subset of $E$ and $T: G\rightarrow G$ a nonexpansive
		mapping. Then $AC(\mathcal{E},G)$ is $T$-invariant (namely,
		$T(AC(\mathcal{E},G))\subset AC(\mathcal{E},G)$), where
		$\mathcal{E}=\{T^m(G), m\in \mathbb{N}\}$.
	\end{lem}
	\begin{proof}
		It is clear that $\mathcal{E}=\{T^m(G), m\in \mathbb{N}\}$ is a decreasing
		sequence of subsets of $G$, then $AC(\mathcal{E},G)$ is a nonempty
		$\mathcal{T}_{\varepsilon,\lambda}$-closed $L^0$-convex subset of $G$ by 
		Lemma \ref{lem.4.4}.
		
		For any $x\in AC(\mathcal{E},G)$, we have
		\begin{align}\nonumber
			AR(\mathcal{E},Tx)&=\bigwedge_{m}R(T^m(G),Tx)\nonumber\\
			&\leq \bigwedge_{m}R(T^{m-1}(G),x)\nonumber\\
			&=AR(\mathcal{E},x)\nonumber\\
			&=AR(\mathcal{E},G),\nonumber
		\end{align}
		which implies that $Tx\in AC(\mathcal{E},G)$. Thus,
		$AC(\mathcal{E},G)$ is $T$-invariant.
	\end{proof}

	\begin{proof}[\textbf{Proof of Theorem \ref{thm.1.16}}]
		Let $\mathcal{E}=\{T^m(G), m\in \mathbb{N}\}$, then
		$AC(\mathcal{E},G)$ is an $L^0$-convexly compact subset with random
		normal structure by 
		Lemma \ref{lem.4.4}. Lemma \ref{lem.4.5} shows that $T$ maps
		$AC(\mathcal{E},G)$ into $AC(\mathcal{E},G)$, further, by
		Theorem \ref{thm.1.15}, there exists $z\in AC(\mathcal{E},G)$ such that
		$Tz=z$. It suffices to show that $z\in C(G,G)$.
		
		By the isometry of  $T$, we have
		$$R(T^m(G),z)=R(T^m(G),Tz)=R(T^{m-1}(G),z), \forall ~m\in \mathbb{N}.$$
		It follows that  $AR(\mathcal{E},G)=AR( \mathcal{E},z )=\bigwedge_{m}
		R(T^m(G),z)=R(G,z)$.
		
		For any $x\in G$, we have
		$$R(G,x)\geq \bigwedge_{m} R(T^m(G),x)=AR( \mathcal{E},x )\geq
		AR(\mathcal{E},G)=R(G,z),$$
		which implies that $R(G,G)\geq R(G,z)$.
		Thus $z\in C(G,G)$.
	\end{proof}
	
    \par As in the proof of Corollary \ref{cor.2.8},  $G$ in Corollary \ref{cor.4.6} is $L^{0}$-convexly compact with random normal structure. Thus,  Corollary \ref{cor.4.6} follows directly from Theorem \ref{thm.1.16}. 
	
	\begin{cor}\label{cor.4.6}
		Let $( E , \| \cdot \|)$ be a $\mathcal{T}_{\varepsilon,\lambda}$-complete
		random uniformly convex $RN$ module over $\mathbb{K}$ with base
		$(\Omega,\mathcal{F},P)$ and $G$ be a nonempty a.s. bounded
		$\mathcal{T}_{\varepsilon,\lambda}$-closed $L^0$-convex subset of $E$. Then
		every isometric mapping $T: G\rightarrow G$ has a fixed point in
		$C(G,G)$.
	\end{cor}

	\par To prove Theorem \ref{thm.1.17}, we first present Lemma \ref{lem.4.7} below.
	
	\begin{lem} \label{lem.4.7}
		Let $( E , \| \cdot \|)$ be a $\mathcal{T}_{\varepsilon,\lambda}$-complete
		$RN$ module over $\mathbb{K}$ with base $(\Omega,\mathcal{F},P)$, $G$
		an $L^0$-convexly compact subset of $E$ and $T: G\rightarrow G$ a surjective isometric mapping. Then $T(C(G,G))=C(G,G)$.
	\end{lem}
	\begin{proof}
		By Proposition \ref{prop.2.3}, $C(G,G)\neq \emptyset$. For any $x\in C(G,G)$, since T is a surjective isometric
		mapping, one has
		$$R(G,Tx)=R(T(G),Tx)=R(G,x)=R(G,G),$$
		and then $Tx\in C(G,G)$. So, $T(C(G,G))\subset
		C(G,G)$.
		
		On the other hand, applying the above argument to $T^{-1}$, we can obtain
		$C(G,G) \subset T(C(G,G)) $. ~
	\end{proof}

	\begin{proof}[\textbf{Proof of Theorem \ref{thm.1.17}}]
		By Proposition \ref{prop.2.3}, $C(G,G)$ is a nonempty $ \mathcal{T}_{\varepsilon,\lambda} $-closed $L^0$-convex subset of $G$. 
		Let $\mathcal{A}=\{ H : H ~$is a nonempty $ \mathcal{T}_{\varepsilon,\lambda} $-closed $L^0$-convex subset of$~ C(G,G) ~$such 
		that$ ~ T(H)\subset H,~ \forall~ T\in \mathcal{T} \}$. Since  $T(C(G,G))=C(G,G)$ 
		for every $T\in \mathcal{T}$ by Lemma \ref{lem.4.7}, then $\mathcal{A}$ is nonempty.
		The $L^0$-convex compactness of $G$ and the Zorn' Lemma together show that there exists a minimal element $H_0$ in $\mathcal{A}$.
		
		We assert that $H_0$ consists of a single element, which is the common fixed 
		point of $\mathcal{T}$. Otherwise, let $A=(D(H_0)>0)$, then $P(A)>0$. Since $G$ has 
		random normal structure, there exists $z_{0}\in H_0$ such that $l:=R(H_{0},z_{0})<D(H_0)$ on $A$. Let
		$$H_1=\{x\in H_0: R(H_0,x)\leq l\},$$
		it is clear that $H_1$ is nonempty, $ \mathcal{T}_{\varepsilon,\lambda} $-closed and $L^0$-convex.
		
		We assert that $H_{1}\in \mathcal{A}$. It suffices to show that $T(H_1)\subset H_1$ 
		for every $ T\in \mathcal{T}$. Observe that $T\in \mathcal{T}$ implies that $T^{-1}\in \mathcal{T}$ 
		and $T^{-1}(H_{0})\subset H_{0}$. For any $T\in \mathcal{T}$ and any $x\in H_1$, we have
		\begin{align}\nonumber
			R(H_0,Tx)&=\bigvee\{\|Tx-h\|:h\in H_0\} \nonumber\\
			&=\bigvee\{\|Tx-T(T^{-1}h)\|:h\in H_0\} \nonumber\\
			&=\bigvee\{\|x-T^{-1}h\|:h\in H_0\} \nonumber\\
			&\leq R(H_0,x) \nonumber\\
			&\leq l. \nonumber
		\end{align}
		Then, $Tx\in H_1$ and hence $T(H_1)\subset H_1$.
		
		By the minimality of $H_0$, we have $H_0=H_1$, and then $D(H_0)\leq l$, which contradicts with the fact that $l<D(H_0)$ on $A$.
	\end{proof}

	\begin{cor}\label{cor.4.8}
		Let $( E , \| \cdot \|)$ be a $\mathcal{T}_{\varepsilon,\lambda}$-complete
		random uniformly convex $RN$ module over $\mathbb{K}$ with base
		$(\Omega,\mathcal{F},P)$, $G$ a nonempty a.s. bounded
		$\mathcal{T}_{\varepsilon,\lambda}$-closed $L^0$-convex subset of $E$ and
		$\mathcal{T}=\{T ~|~ T: G\rightarrow G~$is a surjective isometric mapping$\}$,
		then $\mathcal{T}$ has a common fixed point in $C(G,G)$.
	\end{cor}

	\section{Proofs of Theorems \ref{thm.1.18}, \ref{thm.1.19} and \ref{thm.1.20}}\label{sec.5}

	\begin{proof}[\textbf{Proof of Theorem \ref{thm.1.18}}]
		By Theorems 2.11 and 2.16 of \cite{GZWY20},  $L^0(\mathcal{F},V)$ is an
		$L^0$-convexly compact subset with random normal structure of the
		$\mathcal{T}_{\varepsilon,\lambda}$-complete $RN$ module $L^0(\mathcal{F},B)$.
		
		For each $T\in \mathcal{T}$, define $\hat{T}: L^0(\mathcal{F},V)\rightarrow
		L^0(\mathcal{F},V)$ by
		$$\hat{T}(x)= \text{the equivalence class of}~ T(\cdot, x^0(\cdot)), ~\forall~
		x\in L^0(\mathcal{F},V),$$
		where $x^0$ is an arbitrarily chosen representative of $x$. It is easy to check that $ \hat{T}(x)$ is well defined and $\hat{\mathcal{T}}:=\{\hat{T}: T\in \mathcal{T} \}$ is a commutative family of nonexpansive mappings.
		
		Applying Theorem \ref{thm.1.15} to $\hat{\mathcal{T}}$ and $L^0(\mathcal{F},V)$,
		there exists $x\in L^0(\mathcal{F},V)$ such that $\hat{T}(x)=x$ for each
		$\hat{T} \in \hat{\mathcal{T}}$. Then an arbitrarily chosen representative $x^0$
		of $x$ must satisfy $T(\omega,x^0(\omega) )=x^0(\omega)$ for almost all
		$\omega \in \Omega$ and all $T\in \mathcal{T}$.
	\end{proof}

	\par To prove Theorems  \ref{thm.1.19} and \ref{thm.1.20}, we first give 
	the following Lemma \ref{lem.5.1}. Here, we would like to remind the readers 
	that if $V$ is a nonempty weakly compact convex subset of a Banach space, then $c(V,V)$ 
	is a nonempty closed convex subset of $V$  (see \cite[Proposition 1]{Lim2003} for details).

	\par Conventionally, by identifying each $v\in V$ with $ \tilde{I}_{\Omega}v $, we can regard $V$ as a subset of  $ L^{0}(\mathcal{F},V) $.

	\begin{lem}\label{lem.5.1}
		Let $(B,\|\cdot\|)$ be a Banach space over $\mathbb{K}$ and $V$ be a nonempty weakly compact
		convex subset of $B$. Then 
		$$C(L^{0}(\mathcal{F},V),L^{0}(\mathcal{F},V))=L^{0}(\mathcal{F},c(V,V)).$$
	\end{lem}
	\begin{proof}
		For any $x\in L^{0}(\mathcal{F},V)$, it is  clear that $r(V,x^{0}(\cdot)):\Omega \rightarrow \mathbb{R}_{+}$ is a representative of $R(L^{0}(\mathcal{F},V),x)$, where $x^{0}$ is an arbitrarily chosen representative of $x$. Let $S(\mathcal{F},V)$ be the set of equivalence classes of simple random elements of $ \Omega $ to $ V $,  since $S(\mathcal{F},V)$ is $ \mathcal{T}_{\varepsilon,\lambda}$-dense in $ L^{0}(\mathcal{F},V) $, $R(L^{0}(\mathcal{F},V),\cdot):E\rightarrow L^{0}_{+}(\mathcal{F})$ is $ \mathcal{T}_{\varepsilon,\lambda} $-continuous and $ \|\cdot\|:E\rightarrow L^{0}_{+}(\mathcal{F})$ is also $ \mathcal{T}_{\varepsilon,\lambda} $-continuous, we have
		\begin{align}
			R(L^{0}(\mathcal{F},V),x)&=R(L^{0}(\mathcal{F},V),L^{0}(\mathcal{F},V))\nonumber\\
			&=\bigwedge_{y\in L^{0}(\mathcal{F},V)}R(L^{0}(\mathcal{F},V),y)\nonumber\\
			&=\bigwedge_{\sum_{i=1}^{n}\tilde{I}_{A_{i}}v_{i}\in S(\mathcal{F},V)}R(L^{0}(\mathcal{F},V),\sum_{i=1}^{n}\tilde{I}_{A_{i}}v_{i})\nonumber\\
			&=\bigwedge_{\sum_{i=1}^{n}\tilde{I}_{A_{i}}v_{i}\in S(\mathcal{F},V)}\sum_{i=1}^{n}\tilde{I}_{A_{i}}R(L^{0}(\mathcal{F},V),v_{i})\nonumber\\
			&=\bigwedge_{v\in V}R(L^{0}(\mathcal{F},V),v)\nonumber\\
			&=\bigwedge_{v\in V}\bigvee_{u\in V}\|u-v\|\nonumber\\
			&=R(V,V),\nonumber
		\end{align}
		which, by noting that $ r(V,V) $ is  a representative of $ R(V,V) $,  implies that
		 $$r(V,x^{0}(\cdot))=r(V,V)~ a.s. \text{~iff~} R(L^{0}(\mathcal{F},V),x)=R(L^{0}(\mathcal{F},V),L^{0}(\mathcal{F},V)).$$
		  Hence, the proof is complete.
	\end{proof}

	\begin{proof}[\textbf{Proof of Theorem \ref{thm.1.19}}]
	Define $\hat{T}: L^0(\mathcal{F},V)\rightarrow L^0(\mathcal{F},V)$ by
	$$\hat{T}(x)= \text{the equivalence class of}~ T(\cdot, x^0(\cdot)), ~\forall~
	x\in L^0(\mathcal{F},V),$$
	where $x^0$ is an arbitrarily chosen representative of $x$. It is easy to check that $\hat{T}$ is isometric and  $L^0(\mathcal{F},V)$ is $L^0$-convexly compact with random normal structure.
	
	Applying Theorem \ref{thm.1.16} to $\hat{T}$ and $L^0(\mathcal{F},V)$,
	there exists $x\in C(L^0(\mathcal{F},V)$,\\
	$L^0(\mathcal{F},V))$ such that $\hat{T}(x)=x$. 
	By Lemma \ref{lem.5.1}, $x\in L^{0}(\mathcal{F},c(V,V))$.
	Let $x^0$ be an arbitrarily chosen representative of $x$, then $x^0(\omega)\in c(V,V)$ and $T(\omega,x^0(\omega) )=x^0(\omega)$ for almost all $\omega \in \Omega$.	
	\end{proof}

	\begin{proof}[\textbf{Proof of Theorem \ref{thm.1.20}}]
	Similar to the proof of Theorem \ref{thm.1.19}, by Theorem \ref{thm.1.17}, $\hat{\mathcal{T}}=\{\hat{T}:T\in \mathcal{T}\}$ has a common fixed point $x$ in $C(L^0(\mathcal{F},V),L^0(\mathcal{F},V))$. Then,  for an arbitrarily chosen representative $x^0$ of $x$, $x^0(\omega)\in c(V,V)$ and $T(\omega,x^0(\omega) )=x^0(\omega)$ for almost all $\omega \in \Omega$ and all $T\in \mathcal{T}$.
	\end{proof}

	\section{Concluding remarks}\label{sec.6}
	
	\begin{rem}\label{rem.6.1}
		As in \cite{Guo99,GMT24}, we may suppose that the base space of an $RN$ module $(E,\|\cdot\|)$ is an arbitrary $\sigma$-finite measure space $(\Omega,\mathcal{F},\mu)$. A probability measure $P_{\mu}$ associated with $(\Omega,\mathcal{F},\mu)$ is defined as follows. When $(\Omega,\mathcal{F},\mu)$ is a finite measure space, $P_{\mu}(A)=\frac{\mu(A)}{\mu(\Omega)}$ for each $A\in \mathcal{F}$; when $(\Omega,\mathcal{F},\mu)$ is a $\sigma$-finite measure space, for example, let $\{A_{n},n\in \mathbb{N}\}$ be a countable partition of $\Omega$ to $\mathcal{F}$ such that $0<\mu(A_{n})<+\infty$ for each $n$ in $\mathbb{N}$, then $P_{\mu}(A)=\sum_{n=1}^{\infty}\frac{\mu(A\cap A_{n})}{2^{n}\mu(A_{n})}$ for each $A\in \mathcal{F}$. It is clear that $P_{\mu}$ and $\mu$ are equivalent. When we regarded the $RN$ module $(E,\|\cdot\|)$ as the one with base $(\Omega,\mathcal{F},P_{\mu})$, in this case a sequence $\{x_n, n \in \mathbb{N}\}$ in $(E,\|\cdot\|)$ converges in the $(\varepsilon,\lambda)$-topology to $x$ iff $\{\|x_n -x\|,n \in
		\mathbb{N}\}$ converges in probability measure $P_{\mu}$ to 0, that is, iff $\{\|x_n -x\|,n \in \mathbb{N}\}$ converges locally in measure $ \mu $ to 0. From the equivalence of $P_{\mu}$ to $\mu$, one can easily see that all the results of this paper are still true in the case when the base space is $\sigma$-finite. Further, define $ \interleave \cdot \interleave:E\rightarrow [0,+\infty)$ by $ \interleave x \interleave=\sum_{n=1}^{\infty}\frac{1}{2^{n}\mu(A_{n})}\int_{A_{n}}\|x\|\wedge 1 d\mu$ (or $ \interleave x \interleave=\sum_{n=1}^{\infty}\frac{1}{2^{n}\mu(A_{n})}\int_{A_{n}}\frac{\|x\|}{\|x\|+1} d\mu$) for any $ x\in E $, then it is known that $ (E,\interleave \cdot \interleave) $ is a Fr\'{e}chet space when $(E,\|\cdot\|)$ is $ \mathcal{T}_{\varepsilon,\lambda} $-complete, and the $(\varepsilon,\lambda)$-topology is exactly the linear topology induced by the quasinorm $\interleave \cdot \interleave$. In nonsmooth differential 
		geometry on metric measure spaces, the metric measure spaces are often assumed to be $\sigma$-finite, the $L^{0}$-normed $L^{0}$-modules (namely, $RN$ modules) are just endowed with the linear topology induced by the quasinorm $\interleave \cdot \interleave$ \cite{Gigli2018,LP19,LPV24,CLP2024}. We believe that
		the fixed point theorems developed in \cite{GWXYC24} and this paper will be useful in the future study of nonsmooth differential 
		geometry on metric measure spaces.  
	\end{rem}

	\begin{rem}\label{rem.6.2}
		Geometry of $RN$ modules began with \cite{GZ10,GZ12}, where random uniform and strict convexities were introduced and their relations with classical convexities of the $L^{p}$-spaces generated by the $RN$ modules were established. Recently, these works have been used by Pasqualetto et al. \cite{LPV24} in the study of Banach bundles. Random normal structure and random complete normal structure were studied in \cite{GZWG20} and this paper. But, up to now, the notions of random uniform smoothness and random uniform normal structure for complete $RN$ modules have not even  been defined, we hope to see a good theory of them in the near future since metric fixed point theory in random functional analysis will unavoidably touch on their researches.
	\end{rem}

	\begin{rem}\label{rem.6.3}
		Lau et al. \cite{LZ2008,LZ2012} have made much remarkable contribution to the common fixed point theorems for semigroups of nonexpansive mappings in Banach spaces. In the future, an important research topic is naturally the study of the common fixed point theorems for semigroups of nonexpansive mappings in complete $RN$ modules, but this will be a difficult topic since $RN$ modules are not locally convex in general and this topic will unavoidably involve the stable random weak or weak* topology under the theory of random conjugate spaces, which together with stable compactness is currently the most difficult subject being developed, see \cite{GWT23,GMT24,GWXYC24} for details.	
\end{rem}

	\section*{Acknowledgment}
	
	The first three authors were supported by the National Natural Science Foundation of China
	(Grant No.12371141) and the Natural Science Foundation of Hunan Province of
	China (Grant No.2023JJ30642). The fourth author was supported by the National Natural Science Foundation of China (Grant
	No. U1811461).

\end{document}